\documentclass[twoside,11pt]{article}

\usepackage{jmlr2e}
\usepackage{amsmath,amsthm}
\usepackage[english]{babel}
\usepackage{color}
\usepackage{algorithm2e}  
\usepackage{hyperref}\hypersetup{colorlinks=true,linkcolor=blue} 
\usepackage{multirow}                                                                                     \usepackage{tikz,pgfplots,float}
\pgfplotsset{compat=1.7}
\usetikzlibrary{arrows,backgrounds,positioning,fit}

\def\bfR{{\mathbb R}}
\def\pp{\mbox{\bf P}}
\def\ee{\mbox{\bf E}}
\def\dd{\mbox{\rm d}}
\newcommand*{\Vmax}{ V_{max}\raisebox{+7pt}{\makebox[2ex]{\hspace{-25pt}$\scriptstyle{\bf x}_M$}}\hspace{-9pt} }
\newcommand*{\xM}{ \raisebox{+7pt}{\makebox[2ex]{$\scriptstyle{\bf x}_M$}} }
\newcommand*{\xMone}{ \raisebox{+7pt}{\makebox[4ex]{$\scriptstyle{\bf x}_{M+1}$}} }
\newcommand*{\xMtwo}{ \raisebox{+7pt}{\makebox[4ex]{$\scriptstyle{\bf x}_{M+2}$}} }
\newcommand*{\xMN}{ \hspace{-2pt}\raisebox{+7pt}{\makebox[4ex]{$\scriptstyle{\bf x}_{M+N}$}} }
\newcommand*{\xMNN}{ \hspace{-2pt}\raisebox{+7pt}{\makebox[7ex]{$\scriptstyle{\bf x}_{M+N_0-1}$}} }
\newcommand*{\Upx}{ \hspace{-3pt}\raisebox{+4pt}{\makebox[1ex]{$\scriptstyle{\bf x}$}} }

\numberwithin{equation}{section}
\newtheorem{theorem}{Theorem}[section]

\newtheorem{corollary}[theorem]{Corollary}

\newtheorem{remark}[theorem]{Remark}

\ShortHeadings{}{}
\firstpageno{1}

\begin{document}

\title{Optimal Stopping with Trees: The Details}

\author{\name Sigurd Assing \email s.assing@warwick.ac.uk \\
       \addr Department of Statistics\\
       University of Warwick\\
       Coventry, CV4 7AL, UK
       \AND
       \name Xin Zhi \email xin.zhi@warwick.ac.uk \\
       \addr Department of Statistics\\
       University of Warwick\\
       Coventry, CV4 7AL, UK}

\editor{}

\maketitle

\begin{abstract}
The purpose of this paper is two-fold, first, to review a recent method introduced by
S.\ Becker, P.\ Cheridito, and P.\ Jentzen, for solving high-dimensional optimal stopping problems
using deep Neural Networks, second, to propose an alternative algorithm replacing
Neural Networks by CART-trees which allow for more interpretation of the estimated stopping rules.
We in particular compare the performance of the two algorithms with respect to the
Bermudan max-call benchmark example concluding that the Bermudan max-call may not be suitable 
to serve as a benchmark example for high-dimensional optimal stopping problems. We also show how 
our algorithm can be used to plot stopping boundaries.
\end{abstract}

\begin{keywords}
high-dimensional optimal stopping, sequential data,
CART-trees, American Option, Bermudan Option, deep learning
\end{keywords}
\section{Describing the Problem}\label{problem}
This paper builds upon the work presented in \cite{BCJ2019} whose authors used 
neural networks to directly approximate optimal stopping rules---see the paper
and references therein for the motivation of their approach.

Our interest in this approach comes from applying optimal stopping in the context of 
fast changing sequential data. Before explaining what we mean by {\em fast changing},
let us first describe our framework of optimal stopping.

Assume that ordered data,
\[
{\bf x}_1,\dots,{\bf x}_{M-1},{\bf x}_M,{\bf x}_{M+1},\dots,{\bf x}_{M+N},
\]
is given (or will sequentially be given) 
of which, initially, only ${\bf x}_1,\dots,{\bf x}_{M-1}$ are accessible.
Then, a process is started which makes
${\bf x}_{M+\,n},\,n=0,1,\dots,N$, 
available to a decision maker, one after the other,
for checking a reward which is associated with each of these data points.
The reward can only be taken while being checked.
If the decision maker takes the reward, the process stops and no further
reward can be checked. If they do not take it, the process continues.
Eventually, if none of the rewards associated with 
${\bf x}_{M},\dots,{\bf x}_{M+N-1}$
is taken, the process would also stop and
the decision maker would have to take the reward 
associated with the only remaining data point ${\bf x}_{M+N}$.

Of course,
the decision maker's {\em ultimate goal}, which is to stop the process
by taking the highest reward, can never be achieved for sure.
Instead of going for the highest reward, an objective has to be formulated
which takes into account the uncertainty of missing out on higher rewards
because the process of checking them was stopped.

But how can this uncertainty be taken into account?
When the process starts, 
that is checking the reward associated with ${\bf x}_M$,
the only information available is given by
${\bf x}_1,\dots,{\bf x}_M$ which could be used to infer
possible values of ${\bf x}_{M+1},\dots,{\bf x}_{M+N}$ 
before checking the rewards associated with these data points.
This inference would only be relevant for decision making
if the inferred values can be associated with rewards, too.
Therefore, it is usually assumed that rewards can be determined 
by a family of functions, $u(n,\cdot),\,n=0,\dots,N$,
with each function being
defined on a space which is big enough to contain all values
the data under consideration could take.
In this paper that space is going to be $\bfR^D$,
for some possibly very large $D$.
Adding the parameter $n$ to these functions $u(n,\cdot)$
takes into account that the nature of associated rewards might change 
when sequentially going through the data.

Now, for us, inferring possible values of ${\bf x}_{M+1},\dots,{\bf x}_{M+N}$ 
means generating such values using a generative model which was trained
on ${\bf x}_1,\dots,{\bf x}_M$. 
This way, we could generate many samples,
\[
{\bf x}_{M},{\bf x}_{M+1}^{(k)},\dots,{\bf x}_{M+N}^{(k)},\quad k=1,\dots,K,
\]
of a probability distribution which is hopefully close to the
data generating process behind the yet-unobserved true values
${\bf x}_{M+1},\dots,{\bf x}_{M+N}$. 
Note that each sample can be considered a path starting from
the same data point ${\bf x}_{M}$ which, for notational purposes,
is also denoted by ${\bf x}_{M}^{(k)}$, occasionally.
Since $K$ is supposed to be large,
in what follows, we would identify this ensemble of samples with the distribution
provided by the generative model. The sample mean would therefore play
the role of the expectation operator.

Then, for each $k$, we would be able to find
\[
\mbox{$\hat{\nu}^{(k)}$ such that}\hspace{5pt}
u(\,\hat{\nu}^{(k)},{\bf x}_{M+\,\hat{\nu}^{(k)}}^{(k)}\,)
\,=\,
\max_{0\le n\le N}
u(n,{\bf x}_{M+n}^{(k)})
\]
leading to
\[
\Vmax
\,=\,
\frac{1}{K}\sum_{k=1}^K
u(\,\hat{\nu}^{(k)},{\bf x}_{M+\,\hat{\nu}^{(k)}}^{(k)}\,),
\]
which would be the expected maximal reward,
conditional on what we were able to learn 
from ${\bf x}_1,\dots,{\bf x}_M$
about the yet-unobserved true values ${\bf x}_{M+1},\dots,{\bf x}_{M+N}$.

At this point, we should ask whether $\Vmax$ would be a better 
objective for the decision maker than the ill-posed highest reward.
Unfortunately it is not. 
By general theory of optimal stopping\footnote{Our main reference for general theory
of optimal stopping is \cite{PS2006}, in what follows.},
the optimal stopping rule with respect to the objective $\Vmax$
would still be to stop by taking the highest reward
(as is almost clear from how $\hat{\nu}^{(k)}$ was defined).
So, even in an ideal world of no modelling error\footnote{The generative model
might not capture the full reality.},
it would not be possible to achieve the expected maximum reward
because this would still require knowing all rewards associated with
the yet-unobserved data points ${\bf x}_{M+1},\dots,{\bf x}_{M+N}$
when just checking ${\bf x}_{M}$.

But what type of expected reward could be achieved?
Observe that the objective $\Vmax$ can also be written as
\[
\sup_{\{\nu^{(k)}\}_{k=1}^K}
\frac{1}{K}\sum_{k=1}^K
u(\,{\nu}^{(k)},{\bf x}_{M+\,{\nu}^{(k)}}^{(k)}\,),
\]
where the supremum is over all stopping-ensembles $\{\nu^{(k)}\}_{k=1}^K$ 
the members of which are functions
$\nu^{(k)}:(\bfR^D)^{N+1}\to\{0,\dots,N\}$.
It turns out that, if one puts a constraint 
on the members of possible stopping-ensembles,
it would be possible to achieve the maximal expected reward subject to this constraint.

Define
\[
V_0\xM
\,=\,
\sup_{\{\nu_0^{(k)}\}_{k=1}^K}
\frac{1}{K}\sum_{k=1}^K
u(\,{\nu}_0^{(k)},{\bf x}_{M+\,{\nu}_0^{(k)}}^{(k)}\,),
\]
where the supremum is over all stopping-ensembles $\{\nu_0^{(k)}\}_{k=1}^K$
the members of which can be represented\footnote{
These representations are not unique.}
by 
\[
\nu_0^{(k)}
\,=\,
\sum_{n=1}^N 
n\times
f_{0,n}(\,{\bf x}_{M},{\bf x}_{M+1}^{(k)},\dots,{\bf x}_{M+\,n}^{(k)}\,),
\]
using functions $f_{0,n}:(\bfR^D)^{n+1}\to\{0,1\}$ 
which do \underline{not} dependent on $k$.

Usually, $V_0\xM$ would be smaller than $\Vmax$ but it is as large as it can be
subject to only using the information available at the point of stopping.
Furthermore, general theory of optimal stopping tells that
the supremum used to define $V_0\xM$ is attained
at an optimal stopping-ensemble which can be represented by
a collection of functions ${f}^\star_{0,n},\,n=1,\dots,N$, such that
\[
\sum_{n=1}^N
{f}^\star_{0,n}(\,{\bf x}_{M},{\bf x}_{M+1}^{(k)},\dots,{\bf x}_{M+\,n}^{(k)}\,)
\,=\,
{\bf 1}_{V_0^{{\bf x}_M}\,>\,u(0,{\bf x}_M)}.
\]

Thus, subject to modelling error,
when checking the reward associated with ${\bf x}_{M}$,
the decision
\begin{equation}\label{starRule}
\mbox{\sc do not take this reward {\it iff}\, it is smaller than $V_0\xM$}
\eqno(\star)\nonumber
\end{equation}
would be optimal for achieving the objective given by $V_0\xM$
without the need to know any of the rewards associated 
with yet-unobserved data points.
And if this decision leads to checking the reward associated with ${\bf x}_{M+1}$,
then the decision maker could
\begin{itemize}
\item
re-train their generative model
using ${\bf x}_1,\dots,{\bf x}_M,{\bf x}_{M+1}$;
\item
based on a \underline{new}\footnote{
The observed ${\bf x}_{M+1}$ at the beginning of the path
${\bf x}_{M+1},{\bf x}_{M+2}^{(k)},\dots,{\bf x}_{M+N}^{(k)}$
is supposed to indicate
that these samples are different to the samples generated for $V_0^{{\bf x}_M}$,
slightly abusing notation.}
ensemble of samples
\[
{\bf x}_{M+1},{\bf x}_{M+2}^{(k)},\dots,{\bf x}_{M+N}^{(k)},\quad k=1,\dots,K,
\]
generated by the re-trained generative model, work out
\[
V_1\xMone
\,=\,
\sup_{\{\nu_1^{(k)}\}_{k=1}^K}
\frac{1}{K}\sum_{k=1}^K
u(\,1+{\nu}_1^{(k)},{\bf x}_{M+1\,+\,{\nu}_1^{(k)}}^{(k)}\,)
\]
taking the supremum over all stopping-ensembles $\{\nu_1^{(k)}\}_{k=1}^K$
the members of which can be represented by
\[
\nu_1^{(k)}
\,=\,
\sum_{n=1}^{N-1} 
n\times
f_{1,n}(\,{\bf x}_{M+1},{\bf x}_{M+2}^{(k)},\dots,{\bf x}_{M+1\,+\,n}^{(k)}\,);
\]
\item
make the decision whether to take the reward $u(1,{\bf x}_{M+1})$
by comparing it to $V_1\xMone$\,;
\item
thus move on to checking the next reward or not, and so on, till stopping.
\end{itemize}
\begin{remark}\rm\label{steps}
When the rewards associated with ${\bf x}_{M+\,n},\,n=0,\dots,N$, 
become sequentially available to the decision maker,
if the uncertainty of missing out on higher rewards 
because the process of checking them was stopped
is taken into account by sequentially training a generative model 
for predicting rewards associated with yet-unobserved data points,
then $V_0\xM,V_1\xMone,\dots,V_{N}\xMN$ would be the largest objectives which 
can be achieved by decisions solely based on 
rewards associated with observed data points
and expectations of future rewards.
\end{remark}

Therefore, in the context of machine learning, optimal stopping could come
with two expensive tasks, first, training and re-training a generative model,
second, working out the objectives $V_0\xM,V_1\xMone,\dots,V_{N}\xMN$ in a reliable way.

Re-training the model might be necessary because the data generating process
behind ${\bf x}_{1},\dots,{\bf x}_{M+N}$ has been highly inhomogeneous in the sense
that inferring values for ${\bf x}_{M+\,n+1},\dots,{\bf x}_{M+N}$, 
after the initial training
of the generative model on ${\bf x}_{1},\dots,{\bf x}_{M}$, could follow a distribution
very different to the corresponding predictive distribution of the generative model
after being trained on ${\bf x}_{1},\dots,{\bf x}_{M+\,n}$, 
when further $n\ge 1$ unobserved data points 
have been made available to the decision maker.
Since it is not known at which $n$ the initial training would be out of touch
in the above sense, one would have to re-train the generative model sequentially, 
maybe not at step-size one as we suggested, but definitely at a certain frequency.

If the data generating process behind the given sequential data is likely
to be inhomogeneous as described above, we give these data
the attribute {\em fast changing}.

Assume that the fast changing nature of the data requires the generative model
to be re-trained every $N_0\ll N$ steps. Then, 
omitting the first of the bullet-points above Remark \ref{steps},
working out the objectives $V_0\xM,V_1\xMone,\dots,V_{N}\xMNN$
would be based on $N_0$ ensembles of sample paths starting from
${\bf x}_M,{\bf x}_{M+1},\dots,{\bf x}_{M+N_0-1}$, respectively,
each of them generated by one and the same generative model.

First, generating a full ensemble could come with a significant computational cost,
second, working out practically relevant objectives is usually computationally expensive, 
in particular, when the dimension of the data is large.
So, at least for any period of $N_0$ steps between re-training the generative model,
one would want to estimate decisions for each step right at the beginning of the period
using just one ensemble of sample paths.

That is, when the process of checking rewards starts,
one would want to estimate a stopping rule
associated with the optimal stopping-ensemble at which $V_0\xM$ is attained,
and the decisions given by this estimate would be applied without further adjustments
when checking the next $N_0-1$ rewards. Only then, a new cycle would start
including re-training and re-estimation. 

We therefore conclude that, for our purpose, the problem rather is
to be able to find reliable estimates of optimal stopping rules at low computational cost
than to work out the objectives themselves.
When searching for possible solutions we came across
several recent papers related to this problem.

Most noticeable, we think, is the approach taken in \cite{BCJ2019},
and follow-up papers \cite{BCJ2020,BCJW2021},
where stopping rules are represented by functions similar to our functions
${f}_{0,n},\,n=1,\dots,N$, and these functions are estimated by training deep neural networks.

In \cite{CM2020}, though, stopping rules are modelled by binary decision trees.
However, treating a whole stopping rule as a decision tree comes with high computational cost, 
and the authors prove indeed that estimating such trees is NP-hard, in general.
Nevertheless, an advantage of using decision trees is that the estimated stopping rules 
become interpretable.

In \cite{OCZZ2019}, the authors translate a cost-awareness problem in computing into
an optimal stopping problem and then, similar to \cite{CM2020},
model the corresponding stopping rule by a decision tree.

In \cite{EM2021}, the authors closely follow \cite{BCJ2019} but estimate
the functions ${f}_{0,n},\,n=1,\dots,N$, by successfully applying deep Q-learning.

In \cite{HKRT2021}, the authors improve on the method given in \cite{BCJ2019}
using randomized neural networks though they also discuss the possibility
of applying reinforcement learning in this context.

When experimenting with the algorithm suggested in \cite{BCJ2019},
we noticed that, when implementing it on a standard Desktop Computer,
training the neural networks took particularly long. Furthermore, some
interpretability of the estimated optimal stopping rules would be desirable
in certain applications, see \cite{OCZZ2019} for example.

We therefore came up with the idea of using CART-trees instead of neural networks
for estimating ${f}_{0,n},\,n=1,\dots,N$.
This way, we gain interpretability, on the one hand, 
and we lower computational cost, on the other.
The crucial idea, though, is to estimate each ${f}_{0,n}$ separately,
and this crucial idea goes back to \cite{BCJ2019}, originally, as far as we know.
Instead of modelling the whole stopping rule by a decision tree, as in \cite{CM2020},
our model is built on a forest of trees, which again leads to significantly lower computational cost.

Using CART-trees comes with disadvantages, too, and one of them is the high variance
associated with any tree-method. 
Applying additional variance reduction techniques slows down our algorithm, 
but it remains reasonably faster than the algorithm suggested in \cite{BCJ2019},
mainly because training needs less data,
which is particularly important in the context of fast changing sequential data.

Another problem is dimension: our algorithm becomes inaccurate in certain examples 
involving higher dimensional data. At the moment we do not know whether
this problem is intrinsically tree-related or not. 
There is work in progress to answer this question.

The rest of this paper is organised as follows:
Section \ref{treeli} motivates our main algorithm, that is {\bf Algorithm} \ref{alg1},
but also gives the theoretical background of CART-trees in this context,
while the final Section \ref{appli} deals with applications.

\newpage
\section{The Tree Algorithm}\label{treeli}
It suffices to explain the algorithm for estimating a stopping rule
at which $V_0\xM$ is attained.

We assume that a large number, $K$, of paths
\[
\omega^{(k)}\,=\,
({\bf x}_{M},{\bf x}_{M+1}^{(k)},\dots,{\bf x}_{M+N}^{(k)}),\quad k=1,\dots,K,
\]
have been sampled using a given generative model trained on
${\bf x}_1,\dots,{\bf x}_M$,
and we set
\[
\pp
\,=\,
\sum_{k=1}^K\delta_{\omega^{(k)}}/K
\]
identifying the distribution provided by the generative model
with the ensemble of these samples\footnote{
This measure is also called `empirical measure', later.}.
Here, $\delta_{\omega^{(k)}}$ denotes 
the delta distribution giving mass one to the $k$th path, $\omega^{(k)}$, 
and we consider the probability measure $\pp$ on the sample space 
\[
\Omega\,=\,\{\omega^{(k)}:k=1,\dots,K\}
\]
equipped with the $\sigma$-algebra ${\cal F}$ of all subsets of $\Omega$.
We also define canonical random variables $X_n:\Omega\to\bfR^D,\,n=0,\dots,N$, by
\[
X_0(\omega^{(k)})=\omega^{(k)}_0={\bf x}_M
\quad\mbox{and}\quad
X_n(\omega^{(k)})=\omega^{(k)}_n={\bf x}^{(k)}_{M+\,n},\,n=1,\dots,N,
\]
for each $k=1,\dots,K$.
These random variables generate the natural filtration
\[
{\cal F}_n\,=\,\sigma\{X_{n'}:n'\le n\},\quad n=0,\dots,N,
\]
and by a stopping rule, $\tau$, we mean a random variable 
which maps $\Omega\to\{0,\dots,N\}$ such that
\[
\{\tau=n\}\,\in\,{\cal F}_n,\quad n=0,\dots,N.
\]
Thus, since $\tau=\sum_{n=1}^N n\times{\bf 1}_{\{\tau=n\}}$,
there are functions $f_{0,n}:(\bfR^D)^{n+1}\to\{0,1\}$ such that
\[
\tau\,=\,\sum_{n=1}^N n\times f_{0,n}(X_0,\dots,X_n)
\quad\mbox{and}\quad
{\bf 1}_{\{\tau>0\}}
\,=\,
\sum_{n=1}^N f_{0,n}(X_0,\dots,X_n),
\]
which means that random stopping rules $\tau=\{\tau(\omega^{(k)})\}_{k=1}^K$
are just a different way of writing stopping-ensembles
$\{\nu_0^{(k)}\}_{k=1}^K$ used to define $V_0\xM$ in Section \ref{problem}.

Now, let be given an arbitrary reward function
\[
u:\{0,\dots,N\}\times\bfR^D\to\bfR,
\]
and define
\[
V_{0,n}
\,=\,
\sup_{\tau\in{\cal T}_n}\ee[u(\tau,X_\tau)],
\quad n=0,\dots,N,
\]
where
\[
{\cal T}_n
\,=\,
\{\tau:\mbox{$\tau$ is a stopping rule such that $n\le\tau\le N$}\}.
\]
\begin{remark}\rm
(i)
Because $\Omega$ is a finite set, expectations of all random variables exist.
Moreover, for each $n=0,\dots,N$, the set of stopping rules ${\cal T}_n$ 
is finite, too, and the supremum defining $V_{0,n}$ is actually a maximum.

(ii)
Of course, $V_{0,0}$ coincides with $V_0\xM$ defined in Section \ref{problem},
but $V_{0,1},\dots,V_{0,N}$ should not be confused with $V_1\xMone,\dots,V_N\xMN$.
Indeed, $V_{0,1}$ is calculated using inferred values ${\bf x}_{M+1}^{(k)}$
for a yet-unobserved data point ${\bf x}_{M+1}$, 
while $V_1\xMone$ takes into account the observed data point ${\bf x}_{M+1}$, 
and the values of $V_{0,2}$ \&\ $V_2\xMtwo$ differ in a similar way, and so on.
\end{remark}

However, by general theory of optimal stopping, each $V_{0,n}$ is attained
at an optimal stopping rule $\tau_n\in{\cal T}_n,\,n=0,\dots,N$, where
\[
\tau_N\,=\,N
\]
and
\[
\tau_n
\,=\,
n\times
{\bf 1}_{\{u(n,X_n)\,\ge\,{\bf E}[u(\tau_{n+1},X_{\tau_{n+1}})\,|\,{\cal F}_n]\}}
+
\tau_{n+1}\times
{\bf 1}_{\{u(n,X_n)\,<\,{\bf E}[u(\tau_{n+1},X_{\tau_{n+1}})\,|\,{\cal F}_n]\}},
\]
for $n=0,\dots,N-1$. 

This recursion gives the idea of an algorithm based on
finding the indicator functions
\[
{\bf 1}_{\{u(n,X_n)\,\ge\,{\bf E}[u(\tau_{n+1},X_{\tau_{n+1}})\,|\,{\cal F}_n]\}},
\quad n=0,\dots,N-1,
\]
by working backward from $\tau_N=N$.
But, in general, the $n$th indicator function would still be a function
of $X_0,\dots,X_n$ whose calculation might involve a high numerical cost.
It is therefore desirable to assume an easier situation where
\begin{equation}\label{simplerCondExp}
{\bf E}[u(\tau_{n+1},X_{\tau_{n+1}})\,|\,{\cal F}_n]
\,=\,
{\bf E}[u(\tau_{n+1},X_{\tau_{n+1}})\,|\,{X}_n],
\quad n=0,\dots,N-1,
\end{equation}
which would hold true\footnote{
This follows from combining (1.2.4) and Theorem 1.2 in \cite{PS2006}.}
if $\{X_n\}_{n=0}^N$ was a time-homogeneous Markov chain.
\begin{remark}\rm
In practice, it might be that analysing the output of the generative model
reveals dependencies contradicting the Markov property. 
Then, markovianity could still be achieved by 
combining sample points of a path with `earlier' parts of this path
increasing the dimension $D$. Of course, forcing simulated data to be Markov this way
can only lead to good results if the underlying algorithm performs well in high dimensions.
\end{remark}

For fixed $0\le n\le N-1$, it follows from (\ref{simplerCondExp})
that there is a function $g_n:\bfR^D\to\{0,1\}$ such that
\begin{equation}\label{gn}
{\bf 1}_{\{u(n,X_n)\,\ge\,{\bf E}[u(\tau_{n+1},X_{\tau_{n+1}})\,|\,{\cal F}_n]\}}
\,=\,
g_n(X_n),
\end{equation}
and hence
\[
\tau_n
\,=\,
n\times g_n(X_n)+\tau_{n+1}\times(1-g_n(X_n)).
\]

Assume $\tau_{n+1}$ is known.
Since $V_{0,n}$ is attained at $\tau_n$,
it is also attained at $g_n$ leading to
\begin{equation}\label{sup all g}
V_{0,n}
\,=\,
\sup_{g:\bfR^D\to\{0,1\}}
\ee[\,u(n,X_n)\times g(X_n)+u(\tau_{n+1},X_{\tau_{n+1}})\times(1-g(X_n))\,]
\end{equation}
\[
=\sup_{g:\bfR^D\to\{0,1\}}
\frac{1}{K}\sum_{k=1}^K
[\,u(n,\omega^{(k)}_n)\times g(\omega^{(k)}_n)
+u(\tau_{n+1}(\omega^{(k)}),
\omega^{(k)}_{\tau_{n+1}(\omega^{(k)})})\times(1-g(\omega^{(k)}_n))\,],
\]
where the last line is a reminder of our expectation operator 
being an ensemble average. This supremum over all functions
$g:\bfR^D\to\{0,1\}$ is of course identical to
\begin{equation}\label{max all w}
\max_{{\bf w}}
\frac{1}{K}\sum_{k=1}^{K}
[\,u(n,\omega^{(k)}_n)\times{\rm w}_k
+u(\tau_{n+1}(\omega^{(k)}),
\omega^{(k)}_{\tau_{n+1}(\omega^{(k)})})\times(1-{\rm w}_k)\,],
\end{equation}
where the maximum is over weight-vectors ${\bf w}\in\{0,1\}^{K}$ 
whose components satisfy ${\rm w}_k={\rm w}_{k'}$
for any $1\le k,k'\le K$ such that $\omega^{(k)}_n=\omega^{(k')}_n$.
\begin{remark}\rm\label{max all w discussed}
(i)
It is worth discussing (\ref{max all w}) in the case of $n=0$.
Since $\omega^{(k)}_0={\bf x}_M,\,k=1,\dots,K$,
the maximum would be over two weight-vectors, {\bf w}, only: 
either all components of {\bf w} are zero, or all of them are equal to one.
Thus,
\begin{align*}
V_{0,0}
\,=\,&
\max\{\,u(0,{\bf x}_M)\,,
\frac{1}{K}\sum_{k=1}^K
u(\tau_{1}(\omega^{(k)}),\omega^{(k)}_{\tau_{1}(\omega^{(k)})})\,\}\\
\,=\,&
\max\{\,u(0,{\bf x}_M),V_{0,1}\,\}
\,=\,
\max\{\,u(0,X_0),{\bf E}[u(\tau_{1},X_{\tau_{1}})\,|\,{\cal F}_0]\},
\end{align*}
because ${\cal F}_0$ is the trivial $\sigma$-algebra,
so that (\ref{gn}) implies
\[ 
g_0({\bf x}_M)
\, = \,
{\bf 1}_{u(0,{\bf x}_M)\ge V_{0,1}}
\,=\,
{\bf 1}_{u(0,{\bf x}_M)\ge V_{0,0}}
\,=\,
{\bf 1}
\raisebox{-2pt}
{\makebox[13ex]{
$\scriptstyle{u(0,{\bf x}_M)\,\ge\,V_{0}^{{\bf x}_M}}$
}}.
\]
First, this confirms the decision-rule $(\star)$ stated in Section \ref{problem},
second, the important quantity needed for this decision is $V_{0,1}$.
Note that decision-rule $(\star)$ would be optimal even without assuming that
$(X_n)_{n=0}^N$ is a Markov chain. But without this assumption it is 
computationally harder to work out $V_{0,1}$.
Lastly, the optimal function $g_0$ is highly non-unique:
it is only determined by its value at ${\bf x}_M$ 
and can be freely chosen anywhere else.

(ii)
If $n$ is \underline{not} zero in (\ref{max all w}), and ${\bf w}^\star$ is an argmax,
then ${\rm w}^\star_k,\,k=1,\dots,K$,
would determine the values of an optimal function $g_n$ 
at $\omega^{(k)}_n={\bf x}_{M+\,n}^{(k)}$,
and $g_n$ could again be freely chosen anywhere else, which is
$\bfR^D\setminus\{{\bf x}_{M+\,n}^{(1)},\dots,{\bf x}_{M+\,n}^{(K)}\}$.
However, the difference between $g_0$ and $g_n$ is that 
$g_0$ is determined by a value at an observed data point, which cannot be changed,
while $g_n$ is determined by values at unobserved data points,
which have been sampled using a generative model.
Re-sampling another ensemble, say
$({\bf x}_{M},\tilde{\bf x}_{M+1}^{(k)},\dots,\tilde{\bf x}_{M+N}^{(k)}),\,k=1,\dots,K$,
would ultimately result in another optimal function $\tilde{g}_n$,
which is determined by values at 
$\tilde{\bf x}_{M+\,n}^{(1)},\dots,\tilde{\bf x}_{M+\,n}^{(K)}$,
and which does not give any information on `what to do'
at those points of $\{{\bf x}_{M+\,n}^{(1)},\dots,{\bf x}_{M+\,n}^{(K)}\}$
which were not re-sampled.
\end{remark}

The above remark implies that working out (\ref{max all w})
would not help the decision maker to come up with stopping-predictions
outside the sampled ensemble $\omega^{(k)},\,k=1,\dots,K$,    
which could be considered training data for stopping-predictions.
But, it is clear from (\ref{gn}) that finding a function $g_n({\bf x})$,
which also provides meaningful values at out-of-sample {\bf x},
would give a useful estimate of the function
${\bf 1}
\raisebox{-2pt}{\makebox[10ex]{
$\scriptstyle{u(n,{\bf x})\,\ge\,V_n\Upx}$
}}$ 
needed by the decision maker\footnote{
Apply decision-rule $(\star)$ with respect to $V_n\Upx$.}
when the true value ${\bf x}={\bf x}_{M+\,n}$ is revealed,
in particular when $n\ge 1$.

We therefore suggest 
taking the supremum in (\ref{sup all g}) over a reduced set of functions $g$
with a certain structure. This idea is not new.
In \cite{BCJ2019,BCJ2020,BCJW2021,HKRT2021} for example,
the authors restricted themselves to neural networks, which can still be a large class of possible
functions depending on how complex the neural networks would have to be.
But, if the target is optimal stopping in the context of fast changing sequential data, 
taking the supremum in (\ref{sup all g}) over a class of sufficiently complex neural networks
might be too slow in important applications. 

Our idea is to use a smaller class of functions modelled by CART-trees\footnote{
See Section 9.2 in \cite{HTF2009} for a good introduction to CART-trees.}.
The advantages of choosing such functions are at least two-fold.
First, the computational cost of finding good trees is lower than finding
good neural networks, second, the trees are made of interpretable conditions
while the neural networks would in general be not interpretable.

For our purpose, a CART-tree $g:\bfR^D\to\{0,1\}$ is given by a vector of weights,
${\bf w}\in\{0,1\}^L$, for some positive integer $L$, on the one hand,
and a partition of $\bfR^D$ into hyperrectangles, $R_l,\,l=1,\dots,L$,
on the other. We then represent $g$ by
\[
g({\bf x})
\,=\,
{\rm w}_{q({\bf x})},
\quad\mbox{where}\quad
q({\bf x})
\,=\,
\sum_{l=1}^L l\times{\bf 1}_{R_l}({\bf x}),
\]
and denote the set of all such CART-trees by CART.
\begin{remark}\rm
For given ${\bf w}\in\{0,1\}^L$, and $R_l,\,l=1,\dots,L$,
the following manipulation would not change the corresponding CART-tree:
\begin{itemize}
\item
choose $1\le l\le L$;
\item
cut $R_l$ into two disjoint hyperrectangles $R_{l_1}$ and $R_{l_2}$;
\item
set ${\rm w}_{l_1}={\rm w}_{l_2}={\rm w}_{l}$.
\end{itemize}
\end{remark}

As a consequence, any CART-tree can be turned into a binary decision tree
by applying the above manipulation as often as needed to make 
$R_l,\,l=1,\dots,L$, a recursive binary partition.
For example, Figure 1 shows a partition of $\bfR^2$ into $R_l,\,l=1,\dots,5$,
which is not a recursive binary partition. But, choosing $l=3$ and cutting
$R_l$ into $R_{l_1}$ and $R_{l_2}$ as shown in Figure 2 creates 
a recursive binary partition. 

\begin{figure}[H]
    \centering
    \begin{tikzpicture}[scale = 0.7]
    \begin{axis}[axis background/.style={fill=gray!10},grid=major,
    xmin=-1,xmax=9,ymin=-1,ymax=9,clip = false, axis lines=center,  xtick={1,...,8} , ytick = {1,...,8}, xticklabels={,,3,,5,,},axis line style={draw=none},axis line style={draw=none}, yticklabels={,,3,,5,,}]
    
    \addplot[line width=0.4mm,domain=1:5 , color=blue]{3};
    \addplot[line width=0.4mm,domain=3:7, color=blue]{5};
    \addplot  [line width=0.4mm,mark=none,color = blue] coordinates {(3,3) (3,7)};
    \addplot  [line width=0.4mm,mark=none,color = blue] coordinates {(5,1) (5,5)};
    \node at (axis cs:2,5) {\large $R_4$};
    \node at (axis cs:4,4) {\large $R_5$};
    \node at (axis cs:5,6) {\large $R_3$};
    \node at (axis cs:6,3) {\large $R_2$};
    \node at (axis cs:3,2) {\large $R_1$};
    \addplot[color=black,->] coordinates {(0,0) (8,0)};
    \addplot[color=black,->] coordinates {(0,0) (0,8)};     
    
    \end{axis}
    \end{tikzpicture}
    \hskip 50pt
    \begin{tikzpicture}[scale = 0.7]
    \begin{axis}[axis background/.style={fill=gray!10},grid=major,
    xmin=-1,xmax=9,ymin=-1,ymax=9,clip = false, axis lines=center,  xtick={1,...,8} , ytick = {1,...,8}, xticklabels={,,3,,5,,},axis line style={draw=none},axis line style={draw=none}, yticklabels={,,3,,5,,}]
    
    \addplot[line width=0.4mm,domain=1:5 , color=blue]{3};
    \addplot[line width=0.4mm,domain=3:7, color=blue]{5};
    \addplot[line width=0.4mm,mark=none,color = blue] coordinates {(3,3) (3,7)};
    \addplot  [line width=0.4mm,mark=none,color = blue] coordinates {(5,1) (5,7)};
    \node at (axis cs:2,5) {\large $R_4$};
    \node at (axis cs:4,4) {\large $R_5$};
    \node at (axis cs:4,6) {\large $R_{3_1}$};
    \node at (axis cs:6,6) {\large $R_{3_2}$};
    \node at (axis cs:6,3) {\large $R_2$};
    \node at (axis cs:3,2) {\large $R_1$};
    \addplot[color=black,->] coordinates {(0,0) (8,0)};
    \addplot[color=black,->] coordinates {(0,0) (0,8)};     
    
    \end{axis}
    \end{tikzpicture}\\
    {\bf\tiny  Figure 1}\hspace{5.5cm}{\bf\tiny Figure 2}
\end{figure}

Figure 3 below shows the binary tree
associated with the recursive binary partition shown in Figure 2.
The leaves of this tree `carry' the indices
of the corresponding rectangles. Replacing these indices by the components
of a weight-vector, ${\bf w}\in\{0,1\}^6$, for example ${\bf w}=(0,1,0,0,1,1)^T$,
gives the binary decision tree shown in Figure 4, which is associated with
the CART-tree ${\rm w}_{q({\bf x})}$ where
\[
q({\bf x})
\,=\,
1\times{\bf 1}_{R_1}({\bf x})
+2\times{\bf 1}_{R_2}({\bf x})
+3\times{\bf 1}_{R_{3_1}}({\bf x})
+4\times{\bf 1}_{R_{3_2}}({\bf x})
+5\times{\bf 1}_{R_4}({\bf x})
+6\times{\bf 1}_{R_5}({\bf x}).
\]
\begin{figure}[H]
    \centering
    \begin{tikzpicture}[scale = 0.8,level distance=1cm,level/.style={sibling distance=4cm/#1},background rectangle/.style={fill=gray!10}, show background rectangle]
            \node[blue]{$x[1]\leqslant 5$}
            child[blue] {
                node{$x[2]\leqslant 3$}
                child { node{1}
                }
                child {
                    node{$x[1] \leqslant 3$} 
                    child{node{4}}
                    child{
                        node{$x[2] \leqslant 5$}
                        child{node{5}}
                        child{node{$ 3_1$}}
                        }
                }
                edge from parent
            }
            child[blue] {
                node{$x[2] \leqslant 5$}
                child { node{2} } 
                child { node{$3_2$}}
                edge from parent                   
            };
    \end{tikzpicture}
    \hskip 50pt
    \begin{tikzpicture}[scale = 0.8,level distance=1cm,level/.style={sibling distance=4cm/#1},background rectangle/.style={fill=gray!10}, show background rectangle]
                \node[blue]{$x[1]\leqslant 5$}
                child[blue] {
                    node{$x[2]\leqslant 3$}
                    child { node{0}
                    }
                    child {
                        node{$x[1] \leqslant 3$} 
                        child{node{1}}
                        child{
                            node{$x[2] \leqslant 5$}
                            child{node{1}}
                            child{node{0}}
                            }
                    }
                    edge from parent
                }
                child[blue] {
                    node{$x[2] \leqslant 5$}
                    child { node{1} } 
                    child { node{0}}
                    edge from parent                   
                };
    \end{tikzpicture}\\
{\bf\tiny Figure 3}\hspace{6.5cm}{\bf\tiny Figure 4}
\end{figure}

Mapping a recursive binary partition to a binary decision tree 
is of course not a unique procedure. It starts with the nodes (including the root) 
the choice of which does not have to be unique. 
Furthermore, the conditions at nodes can
be formulated in many ways. In this paper we adopt the following style.
Every node is associated with a component, say $x[d]$, of the input {\bf x},
and a split-point, say $c$,
and we always write the corresponding condition using {\it less than or equal} inequalities,
that is $x[d]\le c$.
Underneath a node, the left branch always refers to the condition being satisfied.
The leaves `carry' the weights associated with the corresponding parts of the partition.
In our case, the weights ONE and ZERO would refer to 
the decisions STOP and CONTINUE, respectively.

Coming back to the original problem of estimating an optimal stopping rule
at which $V_0\xM$ is attained,
we have suggested to evaluate $V_{0,n},\,n=N,\dots,0$, one after the other,
following a backward procedure based on finding stopping rules $\tau_n$
at which each of the $V_{0,n}$ would be attained. Assuming markovianity,
each stopping rule $\tau_n$ could recursively be obtained 
by a function $g_n$ at which the supremum in (\ref{sup all g}) is attained,
so that estimating the optimal stopping rule eventually comes down to
estimating a candidate, $\hat{g}_n$, for each of the functions $g_n,\,n=0,\dots,N-1$.

According to our idea, we replace the right-hand side of (\ref{sup all g}) by
\begin{equation}\label{sup all cart}
\sup_{g\in{\rm CART}}
\ee[\,u(n,X_n)\times g(X_n)+u(\tau_{n+1},X_{\tau_{n+1}})\times(1-g(X_n))\,],
\end{equation}
and approach this problem
by recursively constructing CART-trees $\hat{g}_n,\,n=0,\dots,N-1$,
where $\hat{g}_{N-1}$ has to be constructed first.

Fix the $n$th step, $0\le n\le N-1$, of the recursion.

First, we can assume that stopping rules denoted by 
$\hat{\tau}_{n+1},\dots,\hat{\tau}_{N}$
are known from previous constructions.
If $n=N-1$, we would set $\hat{\tau}_{n+1}=\hat{\tau}_{N}=N$, otherwise
\[
\hat{\tau}_{n+1}
\,=\,
(n+1)\times\hat{g}_{n+1}(X_{n+1})+\hat{\tau}_{n+2}\times(1-\hat{g}_{n+1}(X_{n+1})),
\]
where $\hat{g}_{n+1}$ is of course known from the previous recursion-step.

Plugging $\hat{\tau}_{n+1}$ into (\ref{sup all cart}) yields the objective
\[
\inf_{g\in{\rm CART}}
-\ee[\,u(n,X_n)\times g(X_n)+u(\hat{\tau}_{n+1},X_{\hat{\tau}_{n+1}})\times(1-g(X_n))\,],
\]
where the negative expectation can also be written as
\[
\ee[\,\left(\rule{0pt}{12pt}
u(\hat{\tau}_{n+1},X_{\hat{\tau}_{n+1}})-u(n,X_n)
\right)
\times g(X_n)\,]
-
\ee[\,u(\hat{\tau}_{n+1},X_{\hat{\tau}_{n+1}})\,],
\]
and, since the expectation of $u(\hat{\tau}_{n+1},X_{\hat{\tau}_{n+1}})$ does not depend on $g$,
we in effect arrive at
\begin{equation}\label{cartObjective}
\inf_{g\in{\rm CART}}
\ee[\,\left(\rule{0pt}{12pt}
u(\hat{\tau}_{n+1},X_{\hat{\tau}_{n+1}})-u(n,X_n)
\right)
\times g(X_n)\,].
\end{equation}
\begin{remark}\rm\label{bad cond exp}
When taking the infimum in (\ref{cartObjective}) over neural networks instead of CART-trees,
the optimisation would be over the parameters of the corresponding neural networks.
Finding an optimal CART-tree, though,
requires finding candidates for optimal split-components and split-points,
and these candidates are usually selected from the sample points provided for $X_n$.
We therefore re-write the expectation in (\ref{cartObjective}) as
\[
\ee[\,
\ee\!\left[\rule{0pt}{12pt}
u(\hat{\tau}_{n+1},X_{\hat{\tau}_{n+1}})-u(n,X_n)
\,|\,
X_n\right]
\times g(X_n)\,],
\]
where
\begin{equation}\label{needed}
\ee\!\left[\rule{0pt}{12pt}
u(\hat{\tau}_{n+1},X_{\hat{\tau}_{n+1}})-u(n,X_n)
\,|\,
X_n=\omega^{(k)}_n\right],
\quad k=1,\dots,K,
\end{equation}
can be explicitly calculated, 
since the conditional expectation is with respect to the empirical measure
$\sum_{k=1}^K\delta_{\omega^{(k)}}/K$,
and this is what we are going to do below.
However, it could be questioned whether one should explicitly calculate them as such,
because conditional expectations with respect to empirical measures
can be bad estimates of the corresponding limiting objects.
We will come back to this issue in Section \ref{appli_last}.
\end{remark}

Providing workable values of the conditional expectations in (\ref{needed}),
which are needed for split-finding in this context,
we suggest to perform the below removal-procedure 
at the beginning of any step, $n$, of the backward recursion:
denoting
\begin{equation}\label{defiDelta}
\Delta^{\!(k)}_n
\,=\,
\frac{1}{K}\times
[\,u(\hat{\tau}_{n+1}(\omega^{(k)}),
\omega^{(k)}_{\hat{\tau}_{n+1}(\omega^{(k)})})
-
u(n,\omega^{(k)}_n)\,],          
\quad k=1,\dots,K,
\end{equation}
apply

\vspace{5pt}
\hspace{1.3cm}
\framebox{\parbox[c]{11cm}{
For each $k=1,\dots,K$, find all $k'$ such that $\omega^{(k)}_n=\omega^{(k')}_n$,
associate one $\omega^{(k')}_n$ with $\sum_{k'}\Delta^{\!(k')}_n/\#\{{\rm all}\;k'\}$,
and then remove all other $\omega^{(k')}_n$.}}
\hfill(R)\\

\vspace{5pt}\noindent
to the input $\omega^{(k)}_n,\,\Delta^{\!(k)}_n,\,k=1,\dots,K$.

This removal-procedure\footnote{Its computational cost is comparable to the cost of
providing a proximity matrix for cluster analysis.}
results in a compressed version of the input,
which could be shorter than $K$ and which involves old and new $\Delta$-values
scaled with respect to $K$. We do not need to adapt the scaling to 
the new size of the compressed version because, for split-finding,
only sums of $\Delta$-values will be compared, and not ensemble averages.

When describing split-finding algorithms in what follows,
we are going to employ a notation for input-data
emphasising that the removal-procedure (R) has already been applied:
for $k$ from a subset $I$ of $\{1,\dots,K\}$,
we will use $\pmb{\xi}^{(k)}$ and $\Delta^{\!(k)}$
as substitutes for specific 
sample points $\omega^{(k)}_n$ and corresponding $\Delta$-values, respectively.
The general assumption is therefore that inputs $\pmb{\xi}^{(k)}$
can be considered being different from each other, and that
\[
K\times\Delta^{\!(k)}
\,=\,
\ee\!\left[\rule{0pt}{12pt}                                                                                        u(\hat{\tau}_{n+1},X_{\hat{\tau}_{n+1}})-u(n,X_n)                                                                  \,|\,                                                                                                              X_n=\pmb{\xi}^{(k)}\right],
\quad k\in I.
\]

Next, we present a prototype split-finding algorithm 
repeated applications of which would `grow' a CART-tree
at which the infimum in (\ref{cartObjective}) is attained. 
Being applied at a node,
this algorithm finds a condition describing the split at that node.
The phrase `input-data' refers to data which has been mapped to that node
during previous applications of the algorithm.

Let $\pmb{\xi}^{(k)},\,\Delta^{\!(k)},\,k\in I$, be such input-data,
and let $|I|$ denote the number of elements in the index-set $I$.
If $|I|>1$, the algorithm would carry out two steps,
\begin{itemize}
\item
for each component $d\in\{1,\dots,D\}$, choose a bijection
$\sigma_{\!d}:\{1,\dots,|I|\}\to I$ which makes
$\xi^{(\sigma_{\!d}[k])}[d],\,k=1,\dots,|I|$, an ordered representation of the set
$\{\xi^{(k)}[d]:\,k\in I\}$ satisfying
$\xi^{(\sigma_{\!d}[1])}[d]\le\dots\le\xi^{(\sigma_{\!d}[|I|])}[d]$,
\item
calculate
\[
k_d
\,\in\,
{\rm arg}\max_{1\le k\le|I|-1}
\left(
|\sum_{k'=1}^k\Delta^{(\sigma_{\!d}[k'])}\,|
\vee
|\hspace{-6pt}\sum_{k'=k+1}^{|I|}\Delta^{(\sigma_{\!d}[k'])}\,|
\right),
\quad d=1,\dots,D,
\]
as well as
\[
d^\star
\,\in\,
{\rm arg}\max_{1\le d\le D}
\left(
|\sum_{k'=1}^{k_d}\Delta^{(\sigma_{\!d}[k'])}\,|
\vee
|\hspace{-8pt}\sum_{k'=k_d+1}^{|I|}\Delta^{(\sigma_{\!d}[k'])}\,|
\right),
\]
\end{itemize}
producing four pieces of output: $d^\star,\,k_{d^\star}$, and two index-sets
$I_{left}=\{\sigma_{\!d^\star}[1],\dots,\sigma_{\!d^\star}[k_{d^\star}]\}$
as well as
$I_{right}=\{\sigma_{\!d^\star}[k_{d^\star}+1],\dots,\sigma_{\!d^\star}[|I|]\}$.

This output is used to `grow' two branches at the node,
mapping $\pmb{\xi}^{(k)},\,\Delta^{\!(k)},\,k\in I_{left}$, to the left branch,
and $\pmb{\xi}^{(k)},\,\Delta^{\!(k)},\,k\in I_{right}$, to the right one.
A valid condition describing this split would be
\[
x[d^\star]\le\xi^{(\sigma_{\!d^\star}[k_{d^\star}])}[d^\star].
\]

The new branches could lead to either a new node or to a leaf.
Since the purpose of this simple algorithm is to `grow' an optimal CART-tree
at which the infimum in (\ref{cartObjective}) is attained,
splitting as long as possible might be a successful strategy:
\begin{itemize}
\item
if an index-set has got more than one element,
the corresponding branch would lead to a new node at which the above algorithm is repeated
using the input-data mapped to that branch;
\item
if an index-set has got only one element, say $k$, 
the corresponding branch would lead to a leaf,
which carries weight ZERO if $\Delta^{\!(k)}$ is positive,
and weight ONE otherwise.
\end{itemize}

What would happen when the above algorithm is applied to `grow' a tree?
First, using the removal-procedure (R) 
to process $\omega^{(k)}_n,\,\Delta^{\!(k)}_n,\,k=1,\dots,K$,
yields input-data $\pmb{\xi}^{(k)},\,\Delta^{\!(k)},\,k\in I_{root}$,
at the root-node of the tree\footnote{
The case $|I_{root}|=1$ results in a degenerated tree with only one leaf.}.
Then, repeatedly applying the algorithm would create
a CART-tree, ${\rm w}^\star_{{q}^\star({\bf x})},\,{q}^\star({\bf x})=
\sum_{l=1}^{|I_{root}|} l\times{\bf 1}_{{R}^\star_l}({\bf x})$,
such that, for each $l=1,\dots,|I_{root}|$,
there is exactly one element of $\{\pmb{\xi}^{(k)}:k\in I_{root}\}$ in ${R}^\star_l$,
the index of which is denoted by $k_l$,
and ${\rm w}^\star_l=0$ if $\Delta^{\!(k_l)}$ is positive,
whereas ${\rm w}^\star_l=1$ otherwise.

Obviously, the weighted sum
$\sum_{l=1}^{|I_{root}|}\Delta^{\!(k_l)}\times{\rm w}^\star_l$ 
has got the smallest value of all weighted sums 
$\sum_{k\in I_{root}}\Delta^{\!(k)}\times{\rm w}_k$, with weights ${\bf w}\in\{0,1\}^I$, where
\[
\sum_{l=1}^{|I_{root}|}\Delta^{\!(k_l)}\times{\rm w}^\star_l
\,=\,
\sum_{k\in I_{root}}\Delta^{\!(k)}\times{\rm w}^\star_{{q}^\star(\pmb{\xi}^{(k)})}.
\]
Moreover, by the nature of the removal-procedure (R),
for any function $g:\bfR^D\to\{0,1\}$, it holds that
\[
\sum_{k=1}^K\Delta^{\!(k)}_n\times g(\omega^{(k)}_n)
\,=\,
\sum_{k\in I_{root}}\Delta^{\!(k)}\times g(\pmb{\xi}^{(k)})
\]
leading to
\[
\sum_{k=1}^K\Delta^{\!(k)}_n\times g(\omega^{(k)}_n)
\,=\,
\sum_{k\in I_{root}}\Delta^{\!(k)}\times g(\pmb{\xi}^{(k)})
\,\ge\,
\sum_{k\in I_{root}}\Delta^{\!(k)}\times{\rm w}^\star_{{q}^\star(\pmb{\xi}^{(k)})}
\,=\,
\sum_{k=1}^K\Delta_n^{\!(k)}\times{\rm w}^\star
\hspace{-3pt}\raisebox{-4pt}{\makebox[6ex]{$\scriptstyle{q}^\star(\omega^{(k)}_n)$}}
\hspace{2pt},
\]
that is
\[
\ee[\,\left(\rule{0pt}{12pt} 
u(\hat{\tau}_{n+1},X_{\hat{\tau}_{n+1}})-u(n,X_n)
\right)\times g(X_n)\,]
\,\ge\,
\ee[\,\left(\rule{0pt}{12pt}
u(\hat{\tau}_{n+1},X_{\hat{\tau}_{n+1}})-u(n,X_n)
\right)\times{\rm w}^\star_{{q}^\star(X_n)}\,],
\]
and hence the infimum in (\ref{cartObjective}) is attained at 
the CART-tree ${\rm w}^\star_{{q}^\star(\cdot)}$.

However, since the last inequality holds for all functions $g:\bfR^D\to\{0,1\}$,
the above argument also shows 
that the supremum in (\ref{sup all g}) is attained at that CART-tree, too.
So, the question arises whether the function ${\rm w}^\star_{{q}^\star({\bf x})}$ would be a `good'
candidate for $g_n({\bf x})$ in (\ref{gn}), `good' in the sense of helping the decision maker 
to make stopping predictions outside the sampled ensemble $\omega^{(k)},\,k=1,\dots,K$.

The answer is NO. 
The CART-tree ${\rm w}^\star_{{q}^\star(\cdot)}$ suffers from two drawbacks:
high variance and overfitting.
The combined effect of these two drawbacks is very similar
to what has already been described in Remark \ref{max all w discussed}(ii).

High variance is an intrinsic problem of any tree-based method. 
Reducing variance, though, would not require changing the algorithm itself 
as established techniques are usually based on
taking averages of many trees (bagging) combined with cross-validation.
We will discuss the effectiveness of these techniques further below
when dealing with applications.

But, to avoid overfitting, controlling the size of the tree it `grows' should be part
of the algorithm which, as it stands, rather goes for
maximum sized trees creating a leaf for each single data point.
So, while trying to keep the expectation in (\ref{cartObjective}) as small as possible,
a method which adaptively chooses tree-size from the data should be added.


One purpose of this paper is to demonstrate that, in this context,
a data-driven control of tree-size can be achieved by a minor change of the prototype algorithm,
which is faster and comes at a much lower cost than standard methods like pruning.
In our applications, though, 
we will nevertheless combine this data-driven control with ad hoc controls of node-size
and tree-depth, but results were worse when we solely used ad hoc controls.

Our data-driven control of tree-size only affects the part of the prototype algorithm
where the output is used to `grow' two branches at the node.
Instead, the sizes of
\[
|\sum_{k=1}^{k_{d^\star}}\Delta^{(\sigma_{\!d^\star}[k])}\,|
\vee
|\hspace{-8pt}\sum_{k=k_{d^\star}+1}^{|I|}\Delta^{(\sigma_{\!d^\star}[k])}\,|
\quad\mbox{and}\quad
|\sum_{k\in I}\Delta^{\!(k)}\,|
\]
are compared: if the left term is greater then two branches are grown as in the prototype algorithm,
otherwise the node is turned into a leaf
which carries weight ZERO, if the value of $\sum_{k\in I}\Delta^{\!(k)}$ is positive,
and weight ONE if this value is non-positive.

The new algorithm is shown in {\bf Algorithm} \ref{alg1} below.
Note that the {\bf for}-loop with respect to $k$ allows $\Delta_{left}$ to reach $\Delta$,
and hence $|\Delta|$ is part of the comparison. 
As a consequence, the case $|I|=1$ is covered, too.
\RestyleAlgo{ruled}
\SetKwInput{input}{Input}
\begin{algorithm}\label{alg1}
\small
\DontPrintSemicolon
\caption{$\Delta$-Algorithm for Split Finding}
\input{$I$ (index-set associated with current node)}
\input{$D$ (dimension of data)}
\input{$\pmb{\xi}^{(k)},\,k\in I$ (sampled data points---all different)}
\input{$\Delta^{\!(k)},\,k\in I$ (corresponding increments)}
$\Delta \gets \sum_{k\in I}\Delta^{\!(k)}$\;
$score \gets 0$\;
\For{$d=1$ {\rm to} $D$}{
$\Delta_{left} \gets 0$\;
\For{$k$ {\rm in sorted($I$, by $\xi^{(k)}[d]$)}}{
$\Delta_{left} \gets \Delta_{left}+\Delta^{\!(k)}$\;
$score \gets \max(score,|\Delta_{left}|,|\Delta-\Delta_{left}|)$\;
}}
\eIf{$score>|\Delta|$}{
\rm Split with $x[d^\star]\le\xi^{(k_{d^\star})}[d^\star]$,
where $d^\star$ and $k_{d^\star}$ correspond to max score.
}{
\rm The current node is turned into a leaf 
which carries weight ZERO, if $\Delta>0$, and ONE otherwise.
}
\end{algorithm}

In general, repeated application of the $\Delta$-Algorithm does not `grow' a CART-tree
at which the infimum in (\ref{cartObjective}) is attained. 
For example, at step $n$ of the backward recursion,
consider two-dimensional input-data at the root-node as given in Figure 5.
\begin{figure}[H]
    \centering
    \begin{tikzpicture}[scale = 0.7]
    \begin{axis}[axis background/.style={fill=gray!10},grid=major,
    xmin=-1,xmax=9,ymin=-1,ymax=9,clip = false, axis lines=center,  xtick={1,...,8} , ytick = {1,...,8}, xticklabels={,,3,,5,,},axis line style={draw=none},axis line style={draw=none}, yticklabels={,,3,,5,,}]
       \node[blue,right] at (axis cs:2,6){\normalsize{$\Delta^{(1)} = 2$}};
        \node at (axis cs:2,6) [circle, scale=0.3, draw=blue,fill=blue] {};
        
       \node[blue,right] at (axis cs:3,3){\normalsize{$\Delta^{(3)} =-\frac{1}{2}$}};
        \node at (axis cs:3,3) [circle, scale=0.3, draw=blue,fill=blue] {};
        
       \node[blue,right] at (axis cs:5,5){\normalsize{ $\Delta^{(2)} =-\frac{1}{2}$}};
        \node at (axis cs:5,5) [circle, scale=0.3, draw=blue,fill=blue] {};
        
       \node[blue,right] at (axis cs:6,2){\normalsize{$\Delta^{(4)} = 2$}};
        \node at (axis cs:6,2) [circle, scale=0.3, draw=blue,fill=blue] {};
        \addplot[color=black,->] coordinates {(0,0) (8,0)};
        \addplot[color=black,->] coordinates {(0,0) (0,8)};     
        
    \end{axis}
    \end{tikzpicture}
    \hskip 50pt
    \begin{tikzpicture}[scale = 0.7]
    \begin{axis}[axis background/.style={fill=gray!10},grid=major,
    xmin=-1,xmax=9,ymin=-1,ymax=9,clip = false, axis lines=center,  xtick={1,...,8} , ytick = {1,...,8}, xticklabels={,,3,,5,,},axis line style={draw=none},axis line style={draw=none}, yticklabels={,,3,,5,,}]
       \node[blue,left] at (axis cs:2,6){\normalsize{$\Delta^{(1)} = 2$}};
        \node at (axis cs:2,6) [circle, scale=0.3, draw=blue,fill=blue] {};
        
       \node[blue,right] at (axis cs:3,3){\normalsize{$\Delta^{(3)} = -\frac{1}{2}$}};
        \node at (axis cs:3,3) [circle, scale=0.3, draw=blue,fill=blue] {};
        
       \node[blue,left] at (axis cs:5,5){\normalsize{$\Delta^{(2)} = -\frac{1}{2}$}};
        \node at (axis cs:5,5) [circle, scale=0.3, draw=blue,fill=blue] {};
        
       \node[blue,right] at (axis cs:6,2){\normalsize{$\Delta^{(4)} = 2$}};
        \node at (axis cs:6,2) [circle, scale=0.3, draw=blue,fill=blue] {};
    
        \addplot[line width=0.4mm,domain=2.5:8.5 , color=blue]{2.5};
        \addplot[line width=0.4mm,domain=-0.5:8.5, color=red]{5.5};
        \addplot[line width=0.4mm,mark=none,color = blue] coordinates {(2.5,-0.5) (2.5,8.5)};
        \addplot  [line width=0.4mm,mark=none,color = red] coordinates {(5.5,-0.5) (5.5, 5.5)};
        \addplot[color=black,->] coordinates {(0,0) (8,0)};
        \addplot[color=black,->] coordinates {(0,0) (0,8)};

    \end{axis}
    \end{tikzpicture}\\
    {\bf\tiny Figure 5}\hspace{5.5cm}{\bf\tiny Figure 6}
\end{figure}

Since $|\Delta^{(1)}+\Delta^{(2)}+\Delta^{(3)}+\Delta^{(4)}|$ has got the largest absolute value
of any sum one could form with these $\Delta$-values,
the $\Delta$-Algorithm turns the root into a leaf carrying weight ZERO,
because $\Delta^{(1)}+\Delta^{(2)}+\Delta^{(3)}+\Delta^{(4)}=3$ is positive,
resulting in a degenerated tree: all points in $\bfR^2$ are mapped to weight ZERO.

In contrast, Figure 6 shows two binary partitions, each of which is associated with
a CART-tree at which the infimum in (\ref{cartObjective}) would be attained,
but none of them is found by the $\Delta$-Algorithm. 
Moreover,
the CART-tree associated with the blue partition would map
any ${\bf x}\in\bfR^2$ satisfying $x_1\ge 3$ and $x_2\ge 3$
to weight ONE, 
while any ${\bf x}\in\bfR^2$ with $x_1\le 2$ or $x_2\le 2$
would be mapped to weight ZERO.
The CART-tree associated with the red partition, though, would map  
any ${\bf x}\in\bfR^2$ satisfying $x_1\le 5$ and $x_2\le 5$
to weight ONE, 
while any ${\bf x}\in\bfR^2$ with $x_1\ge 6$ or $x_2\ge 6$
would be mapped to weight ZERO.
Thus, despite both of them being candidates for $g_n$ in (\ref{gn}),
the functions these two trees represent massively differ from each other.

This example gives a good idea on how tree-size is controlled by the $\Delta$-Algorithm 
and how the output would change accordingly:
the pocket of relatively small $\Delta$-values inside the square shown at the 
centre of Figure 6 is ignored, and the whole area spanned by the data points in this pocket
is turned into a part of the continuation region.

On the other hand, if all $\Delta^{\!(k)}$
of input-data $\pmb{\xi}^{(k)},\,\Delta^{\!(k)},\,k\in I$, at a node have the same sign
then the $\Delta$-Algorithm would turn this node into a leaf giving the same weight
to the whole area spanned by the data points $\pmb{\xi}^{(k)},\,k\in I$, 
without ignoring anything, 
while the prototype algorithm would `grow' a much bigger tree creating a leaf for
each single data point.

On the whole, one would hope that those parts of the distribution of the sampled ensemble
$\omega^{(k)},\,k=1,\dots,K$, 
possibly ignored by the $\Delta$-Algorithm are mostly\footnote{
Subject to constructed but uninteresting counterexamples.}
unimportant for decision making, which is hard to prove in theory. 
However, in the next section we will see that
estimating $g_n$ by CART-trees which are built on output of the $\Delta$-Algorithm
is quite effective when being applied in interesting examples.

\newpage
\section{Applications}\label{appli}
Following the approach described in Sections \ref{problem} \&\ \ref{treeli},
assume that a generative model fit to given data
${\bf x}_1,\dots,{\bf x}_M$ 
has been used to generate an ensemble of sample paths
\begin{equation}\label{traini}
\omega^{(k)}\,=\,
({\bf x}_{M},{\bf x}_{M+1}^{(k)},\dots,{\bf x}_{M+N}^{(k)}),\quad k=1,\dots,K,
\end{equation}
for estimating a stopping rule, $\hat{\tau}_0$,
at which $V_0\xM=V_{0,0}$ is attained.

We have suggested working backward from $\hat{\tau}_N=\tau_N=N$ via
\begin{equation}\label{forTauZero}
\hat{\tau}_n\,=\,n\times\hat{g}_n(X_n)+\hat{\tau}_{n+1}\times(1-\hat{g}_n(X_n)),
\quad n=N-1,\dots,0,
\end{equation}
where the functions $\hat{g}_n$ are our estimates of the objects $g_n$ 
associated with (\ref{gn}). Each estimate is  based on our $\Delta$-algorithm
combined with bagging and cross validation to counter the high variance of tree methods.
It follows a short description of how bagging and cross validation 
has been applied in our numerical examples.

First, let $\bigcup_{\beta=1}^B\mathfrak{K}_\beta$ be a uniform decomposition
of the index-set $\{1,\dots,K\}$ splitting the above ensemble of sample paths into 
$B$ equally sized bags.

Then, for all bags, we would `grow' CART-trees, $\hat{g}_{n,\beta}$,
and simply set our estimates to be
\begin{equation}\label{ourEstimates}
\hat{g}_{n}
\,=\,
\chi(\,\frac{1}{B}\sum_{\beta=1}^B\hat{g}_{n,\beta}\,),
\quad n=0,\dots,N-1,
\end{equation}
where $\chi\,=\,{\bf 1}_{[\frac{1}{2},\infty)}$ acts as a projector 
mapping averages to either ZERO or ONE.

But, how would we grow the CART-trees $\hat{g}_{n,\beta}$?

Consider the initial case $n=N-1$.
Since $\hat{\tau}_N=\tau_N=N$ is the same for all $\beta$, the values of
$\Delta^{\!(k)}_{N-1},\,k\in\mathfrak{K}_\beta$, calculated according to (\ref{defiDelta}),
would not depend on $\beta$, either. For each $\beta=1,\dots,B$,
we can therefore apply the removal-procedure (R) to 
$\omega^{(k)}_{N-1},\,\Delta^{\!(k)}_{N-1},\,k\in\mathfrak{K}_\beta$, to obtain
initial input-data at the root-node of a CART-tree $\hat{g}_{N-1,\beta}$,
which is subsequently `grown' to full size by repeatedly applying the $\Delta$-algorithm.

In any of the subsequent cases $n=0,\dots,N-2$, though, we would calculate
\[
\Delta^{\!(k)}_{n,\beta}
\,=\,
\frac{1}{K}\times  
[\,u(\hat{\tau}_{n+1,\beta}(\omega^{(k)}),
\omega^{(k)}_{\hat{\tau}_{n+1,\beta}(\omega^{(k)})})                                                      -                             
u(n,\omega^{(k)}_{n})\,],     
\quad k\in\mathfrak{K}_\beta,
\]
based on\footnote{Setting $\hat{\tau}_{N,\beta}=N$, for all $\beta=1,\dots,B$.}
\begin{align*}
\hat{\tau}_{n+1,\beta}(\omega^{(k)})
&\,=\,
(n+1)\times\chi\left(\rule{0pt}{12pt}\right.
\frac{1}{B-1}\sum_{\beta'\not=\beta}\hat{g}_{n+1,\beta'}(\omega^{(k)}_{n+1})
\left.\rule{0pt}{12pt}\right)\\
&\,+\,
\hat{\tau}_{n+2,\beta}(\omega^{(k)})\times(1-
\chi\left(\rule{0pt}{12pt}\right.   
\frac{1}{B-1}\sum_{\beta'\not=\beta}\hat{g}_{n+1,\beta'}(\omega^{(k)}_{n+1})    
\left.\rule{0pt}{12pt}\right)
),
\end{align*}
where applying $\chi$ to a reduced average, 
that is leaving out $\hat{g}_{n+1,\beta}$ when calculating the $\beta$-bag estimate, 
is how we apply the idea of cross validation in this context.

Once all estimates $\hat{g}_n$ have been calculated according to (\ref{ourEstimates}),
they are recursively plugged into (\ref{forTauZero}), 
which eventually yields\footnote{Setting $\hat{g}_N\equiv N$.}
\begin{equation}\label{ourTauZero}
\hat{\tau}_0
\,=\,
\sum_{n=1}^N n\times\underbrace{
(1-\hat{g}_{0}(X_{0}))\times\cdots\times(1-\hat{g}_{n-1}(X_{n-1}))\times\hat{g}_n(X_n)
}_{
\hat{f}_{0,n}(X_0,\dots,X_n)},
\end{equation}
that is our estimate of a stopping rule at which $V_0\xM=V_{0,0}$ is attained.

Recall that
\begin{equation}\label{trueValue}
V_{0,0}
\,=\,
\mbox{$\sup_{\tau}$}\;\ee[u(\tau,X_\tau)],
\end{equation}
where the supremum is over all stopping rules bounded by $N$,
and the expectation operator coincides with the ensemble average, that is
\[
\ee[u(\tau,X_\tau)]
\,=\,
\frac{1}{K}\sum_{k=1}^K
u(\tau(\omega^{(k)}),\omega^{(k)}_{\tau(\omega^{(k)})}).
\]

When replacing $\tau$ on the above right-hand side by $\hat{\tau}_0$, taking into account
(\ref{ourTauZero}) and the structure of the sample paths given by (\ref{traini}),
we arrive at
\[
V^{train}_{0,0}
\,=\,
u(0,{\bf x}_M)\hat{g}_0({\bf x}_M)
+
\sum_{n=1}^N\frac{1}{K}
\sum_{k=1}^K u(n,{\bf x}^{(k)}_{M+\,n})
\hat{f}_{0,n}({\bf x}_M,{\bf x}^{(k)}_{M+1},\dots,{\bf x}^{(k)}_{M+\,n}).
\]

Here, the superscript `train' is used to emphasise that both
the construction of $\hat{\tau}_0$ and the expectation operator in (\ref{trueValue})
are based on the same ensemble of sample paths also called {\it training data}.
Note that $V_{0,0}$ in (\ref{trueValue}) is the true objective with respect to
training data, while $V^{train}_{0,0}$ is a lower bound because $\hat{\tau}_0$
is only one of the stopping rules the supremum in (\ref{trueValue}) is taken over.

Next, let 
$({\bf x}_{M},\tilde{\bf x}_{M+1}^{(k)},\dots,\tilde{\bf x}_{M+N}^{(k)}),\,k=1,\dots,\tilde{K}$,
be another ensemble of sample paths,
also called {\it test data}, which gives another ensemble average, say $\tilde{\ee}$, leading to
\[
\tilde{V}_{0,0} 
\,=\,
\mbox{$\sup_{\tilde{\tau}}$}\;\tilde{\ee}[u(\tilde{\tau},\tilde{X}_{\tilde{\tau}})],
\]
which is the true objective with respect to test data.
Plugging the random variables $\tilde{X}_n,\,n=0,\dots,N$, into the right-hand side
of (\ref{ourTauZero}) gives one of the stopping rules $\tilde{\tau}$
the above supremum is over, and hence 
\begin{equation}\label{ourLowerBound}      
V^{test}_{0,0}
\,=\,                    
u(0,{\bf x}_M)\hat{g}_0({\bf x}_M)    
+
\sum_{n=1}^N\frac{1}{\tilde{K}}    
\sum_{k=1}^{\tilde{K}}u(n,\tilde{\bf x}^{(k)}_{M+\,n})       
\hat{f}_{0,n}({\bf x}_M,\tilde{\bf x}^{(k)}_{M+1},\dots,\tilde{\bf x}^{(k)}_{M+\,n})  
\end{equation}
is a lower bound of $\tilde{V}_{0,0}$.
\begin{remark}\rm\label{overfitting}
The decision-rule $(\star)$ on page \pageref{starRule} could be considered
a statistical model for optimal stopping when checking the reward associated with ${\bf x}_M$.
This model has got one parameter, which we estimated by $V_0\xM$ 
at the point of formulating the rule in Section \ref{problem}.
When using this estimate, the model of course overfits the training data
because $V_0\xM=V_{0,0}$ is based on taking ensemble averages with respect to training data.
By the same reason, when using $\tilde{V}_{0,0}$ instead, the model would overfit the test data.
\end{remark}

In light of the above remark, applying $(\star)$ with $V_0\xM$ replaced by $V^{test}_{0,0}$
should eventually result in a model which is  much less likely to overfit the test data.
Also, $V^{test}_{0,0}$ seems to be a good quantity for judging the quality of 
our stopping rule $\hat{\tau}_0$. Indeed, if
\[
{\tau}_0^\star
\,=\,
\sum_{n=1}^N n\times{f}_{0,n}^\star(X_0,\dots,X_n)
\]
is another estimate of a stopping rule at which $V_0\xM=V_{0,0}$ is attained,
whose defining functions $f_{0,n}^\star$ were estimated 
using some possibly different training data, then\footnote{
Of course, $g_0^\star\,=\,1-\sum_{n=1}^N f_{0,n}^\star$.}
\begin{equation}\label{theirLowerBound}
V^{\star}_{0,0}
\,=\,
u(0,{\bf x}_M){g}_0^\star({\bf x}_M)
+
\sum_{n=1}^N\frac{1}{\tilde{K}}
\sum_{k=1}^{\tilde{K}}u(n,\tilde{\bf x}^{(k)}_{M+\,n})
{f}_{0,n}^\star({\bf x}_M,\tilde{\bf x}^{(k)}_{M+1},\dots,\tilde{\bf x}^{(k)}_{M+\,n})
\end{equation}
is also a lower bound of the true objective $\tilde{V}_{0,0}$ 
with respect to the given test data.
Since both lower bounds ${V}_{0,0}^\star$ and $V^{test}_{0,0}$ 
are obtained by plugging unseen data into estimated objects, 
the higher of the two lower bounds should be associated with the better estimate 
because the corresponding object adapts better to unseen test data.
\begin{corollary}\label{better estimate}
Observing $V^{test}_{0,0}>V_{0,0}^\star$ is an indication for both:
\begin{itemize}
\item
$\hat{\tau}_0$ being the better estimate than $\tau_0^\star$;
\item
$\hat{f}_{0,n}$ acting better on unseen data than $f^\star_{0,n}$.
\end{itemize}
\end{corollary}

As a matter of fact, only observing $V^{test}_{0,0}>V_{0,0}^\star$ would not give a strong
indication for the second item to be true for all $n=1,\dots,N$,
the indication is rather that it should be true for at least some $n\in\{1,\dots,N\}$.
One could nevertheless introduce and compare ${V}_{0,n}^\star$ and $V^{test}_{0,n}$, for $n=1,\dots,N$,
to harden this statement. In this paper, though, 
we restrict ourselves to comparing ${V}_{0,0}^\star$ and $V^{test}_{0,0}$, only.

In what follows, we analyse output of the $\Delta$-algorithm in three scenarios.

First, in Section \ref{theBoundary},  we show how our estimate $\hat{\tau}_0$ can be used
to plot an approximation of the stopping boundary associated with finite horizon
American puts, which gives insight into the quality of the estimate $\hat{\tau}_0$ itself.

Second, in Section \ref{maxCall1}, we compare our results with the results obtained
in \cite{BCJ2019} where the authors applied their method to pricing 
high-dimensional Bermudan max-call options.

Third, in Section \ref{maxCall2}, we compare our results with the results obtained
in \cite{CM2020} where the authors applied their method to pricing Bermudan max-call options
with barrier.
\subsection{Stopping Boundary of American Puts}\label{theBoundary}
The optimal stopping problem associated with American puts (see \cite{PS2006}, Section 25)
is about analysing (\ref{trueValue}) in continuous time, where the reward function would read
\begin{equation}\label{ctReward}
\exp\{-rt\}(C-x)^+,\quad t\in[0,T],\,x\in\bfR,
\end{equation}
and the random variables $X_t$ would follow a geometric Brownian motion 
solving the stochastic differential equation (SDE)
\[
\dd X_t\,=\,\mu X_t\,\dd t + \sigma X_t\,\dd W_t,\quad X_0={\rm x}_M>0.
\]

Here, $r>0$ is the interest rate, $C>0$ is the strike level, $T>0$ is the maturity, 
$\mu$ is the drift coefficient,
$\sigma>0$ is the volatility coefficient, and $W_t$ stands for a one-dimensional standard
Wiener process. The SDE can be considered the generative model whose coefficients
have been estimated using one-dimensional observed data ${\rm x}_1,\dots,{\rm x}_M$.

For numerical purposes, the problem has to be discretised which we achieve by splitting
the maturity $T$ into $N$ equally long time periods with end points $nT/N,\,n=0,\dots,N$,
leading to one-dimensional training data (\ref{traini}) being generated by
\begin{equation}\label{dynamics}
{\rm x}_{M+\,n}^{(k)}
\,=\,
{\rm x}_M\times\exp\{\,
(\mu-\frac{\sigma^2}{2})\frac{nT}{N}
+\sigma\sqrt{\frac{T}{N}}\sum_{n'=1}^n\varepsilon_{n'}^{(k)}
\,\}
\end{equation}
using i.i.d.\ standard normals $\varepsilon_{n'}^{(k)},\,n'=1,\dots,N,\,k=1,\dots,K$.
Using another but independent batch of i.i.d.\ standard normals,
$\tilde{\varepsilon}_{n'}^{(k)},\,n'=1,\dots,N,\,k=1,\dots,\tilde{K}$,
the test data is going to be generated following the same dynamics.
Eventually, at discretised level, the continuous time reward function (\ref{ctReward})
translates into
\[
u(n,x)\,=\,
\exp\{-rnT/N\}(C-x)^+,\quad n\in\{0,\dots,N\},\,x\in\bfR.
\]

All in all, we are now in the position to work out $V_{0,0}^{train}$ and $V_{0,0}^{test}$
for different choices of $\sigma,\,{\rm x}_M,\,T$, and $N$, using the $\Delta$-algorithm
combined with bagging and cross validation as described at the beginning of Section \ref{appli}.
Interest rate $r=0.05$, drift $\mu=0.05$,
and strike level $C=100$, will be fixed in this numerical experiment.
We generate $K=200000$ training samples, and $\tilde{K}=200000$ test samples.
The number of bags is chosen to be $B=10$,
and we also put a constraint on both  maximum tree-depth $(=10)$
and minimum node-size $(=10)$.

Our results are shown in Table \ref{valueAmerPut} below. 
The values of $V_{0,0}^{train}$ and $V_{0,0}^{test}$ are compared with further values
denoted by LS (for Longstaff-Schwartz), referring to a well-established algorithm---see \cite{LS2001}. 
The Longstaff-Schwartz algorithm is not about estimating
an optimal stopping rule, it instead directly estimates the objective using a method
based on regression with respect to special basis functions.
This algorithm works very well in low dimensions, and hence the values it suggests
can be considered benchmarks for our values in this example.

First, our values match the benchmark-values quite well: 
the average relative error $|V_{0,0}^{train}-LS^{train}|/LS^{train}$
and $|V_{0,0}^{test}-LS^{test}|/LS^{test}$ is $0.16\%$ and $0.22\%$, respectively.
Also, the average relative error
$|V_{0,0}^{train}-V_{0,0}^{test}|/V_{0,0}^{test}$ is only $0.14\%$,
so $V_{0,0}^{train}$ and $V_{0,0}^{test}$ do not appear to be very different.

Interestingly, except for two, 
all $V_{0,0}^{test}$-values are bigger than the $V_{0,0}^{train}$-values.
Why is this relation unexpected?
Note that $V_{0,0}^{train}$ would only be equal to the true objective
with respect to training data, i.e.\ $V_{0,0}$, if the supremum in (\ref{trueValue})
is attained at $\hat{\tau}_0$ given by (\ref{ourTauZero}), which is very unlikely
because finding $\hat{\tau}_0$ involves bagging and cross validation. Nevertheless, 
like\footnote{See Remark \ref{overfitting}.} $V_{0,0}$, 
the values of $V_{0,0}^{train}$ would be prone to overfitting
as they entirely depend on the training data. One would therefore expect that
the lower bound $V_{0,0}^{train}$ comes closer to the true objective with respect to training data,
i.e.\ $V_{0,0}$,
than the lower bound $V_{0,0}^{test}$ would come to the true objective with respect to test data, 
i.e.\ $\tilde{V}_{0,0}$.
Furthermore, since $K=\tilde{K}=200000$ is quite large, 
one would also expect $V_{0,0}$ and $\tilde{V}_{0,0}$ to be close,
so that $V_{0,0}^{train}$ should mostly overshoot $V_{0,0}^{test}$.

\begin{table}[ht]
\centering
\begin{tabular}{cccc|cc|cc}
\hline
         &       &     &     & \multicolumn{2}{c|}{training data} & \multicolumn{2}{c}{test data}  \\
$\sigma$&${\rm x}_M$&$T$&$N$&$V_{0,0}^{train}$&$LS^{train}$&$V_{0,0}^{test}$&$LS^{test}$\\ \hline
0.2      & 85    & 1   & 50  & 15.304             & 15.291        & 15.285           & 15.268      \\
0.2      & 90    & 1   & 50  & 11.477             & 11.449        & 11.482           & 11.458      \\
0.2      & 100   & 1   & 50  & 6.058              & 6.046         & 6.068            & 6.049       \\
0.2      & 110   & 1   & 50  & 2.968              & 2.963         & 2.979            & 2.968       \\
0.4      & 85    & 1   & 50  & 20.819             & 20.808        & 20.849           & 20.823      \\
0.4      & 90    & 1   & 50  & 18.104             & 18.091        & 18.129           & 18.101      \\
0.4      & 100   & 1   & 50  & 13.603             & 13.605        & 13.618           & 13.609      \\
0.4      & 110   & 1   & 50  & 10.168             & 10.160        & 10.172           & 10.166      \\
0.2      & 85    & 1   & 100 & 15.282             & 15.280        & 15.301           & 15.279      \\
0.2      & 90    & 1   & 100 & 11.493             & 11.469        & 11.476           & 11.442      \\
0.2      & 100   & 1   & 100 & 6.076              & 6.056         & 6.081            & 6.063       \\
0.2      & 110   & 1   & 100 & 2.965              & 2.955         & 2.994            & 2.982       \\
0.2      & 85    & 2   & 50  & 15.909             & 15.889        & 15.918           & 15.871      \\
0.2      & 90    & 2   & 50  & 12.565             & 12.532        & 12.573           & 12.528      \\
0.2      & 100   & 2   & 50  & 7.677              & 7.655         & 7.685            & 7.663       \\
0.2      & 110   & 2   & 50  & 4.605              & 4.592         & 4.615            & 4.599       \\
0.4      & 85    & 2   & 50  & 24.256             & 24.240        & 24.283           & 24.251      \\
0.4      & 90    & 2   & 50  & 21.901             & 21.896        & 21.943           & 21.913      \\
0.4      & 100   & 2   & 50  & 17.920             & 17.895        & 17.949           & 17.917      \\
0.4      & 110   & 2   & 50  & 14.687             & 14.666        & 14.692           & 14.681      \\ 
\hline
\end{tabular}
\caption{Comparing $V_{0,0}^{train}$ and $V_{0,0}^{test}$ with each other,
and with their Longstaff-Schwartz counterparts, in the American put example.}
\label{valueAmerPut}
\end{table}

However, it might be wrong to assume consistency, that is expecting
$V_{0,0}$ and $\tilde{V}_{0,0}$ to be close to each other.
They both do depend on expectations with respect to different empirical measures,
and these empirical measures would converge weakly to a limit when $K,\tilde{K}\to\infty$.
But the functionals defining $V_{0,0}$ and $\tilde{V}_{0,0}$ do not seem to be continuous enough,
and hence weak convergence of measures might not be sufficient to justify consistency.

Another observation is that all $V_{0,0}^{test}$-values are bigger than the $LS^{test}$-values.
Again, the larger $V_{0,0}^{test}$ 
the better the estimated functionals defining $\hat{\tau}_0$ adapt to unseen data.
On the other hand, 
since the regression function defining $LS^{test}$ has been applied to the same unseen data,
similar to the reasoning leading to Corollary \ref{better estimate}, 
we could say that our estimate $V_{0,0}^{test}$ performs better than $LS^{test}$.

Next, we turn to analysing $\hat{\tau}_0$ itself. 
It is well-known that the continuous time optimal stopping rule associated with American puts
can be described by the underlying geometric Brownian motion, $X_t$,
hitting a stopping boundary, $b(t)$, 
that is the stopping rule can be written as
\[
\inf\{0\le t\le T:X_t\le b(t)\},
\]
where $b:[0,T]\to\bfR$ is a function.
The function $b(t)$ can be found by solving a non-linear Volterra equation of the second kind,
see \cite[Section 25.2]{PS2006} for further details. For this paper,
we solved the Volterra equation by applying the numerical method presented in \cite{K2008}.

Figure 7a shows the solution for $r=\mu=0.05,\,C=100,\,T=1,\,\sigma=0.2,\,{\rm x}_M=85$.
We call it {\it theoretical boundary} because it refers to the continuous time stopping rule,
and by $b(n)$ we mean the continuous time solution $b(t)$ 
being evaluated at discrete time points $t=nT/N,\,n=0,\dots,N$, using $N=50$ for the plot.

We then generate training data
$\omega^{(k)}=({\rm x}_{M},{\rm x}_{M+1}^{(k)},\dots,{\rm x}_{M+N}^{(k)})$,
$k=1,\dots,{K}=200000$,
and test data
$\tilde{\omega}^{(k)}=({\rm x}_{M},\tilde{\rm x}_{M+1}^{(k)},\dots,\tilde{\rm x}_{M+N}^{(k)})$,
$k=1,\dots,\tilde{K}=200000$,
according to (\ref{dynamics}), which yields canonical random variables
$X_0,\dots,X_N$, and $\tilde{X}_0,\dots,\tilde{X}_N$, respectively.
The bullet points underneath the theoretical boundary in Figure 7a refer to 
stopping test sample paths with respect to the stopping rule\footnote{
Bullet points at $N=50$ are not shown as they could be above the strike level $C=100$.}
\begin{equation}\label{notOptimal}
\inf\{0\le n\le N:\tilde{X}_n\le b(n)\}.
\end{equation}
\begin{center} 
\includegraphics[width=0.46\textwidth,height = 5cm]{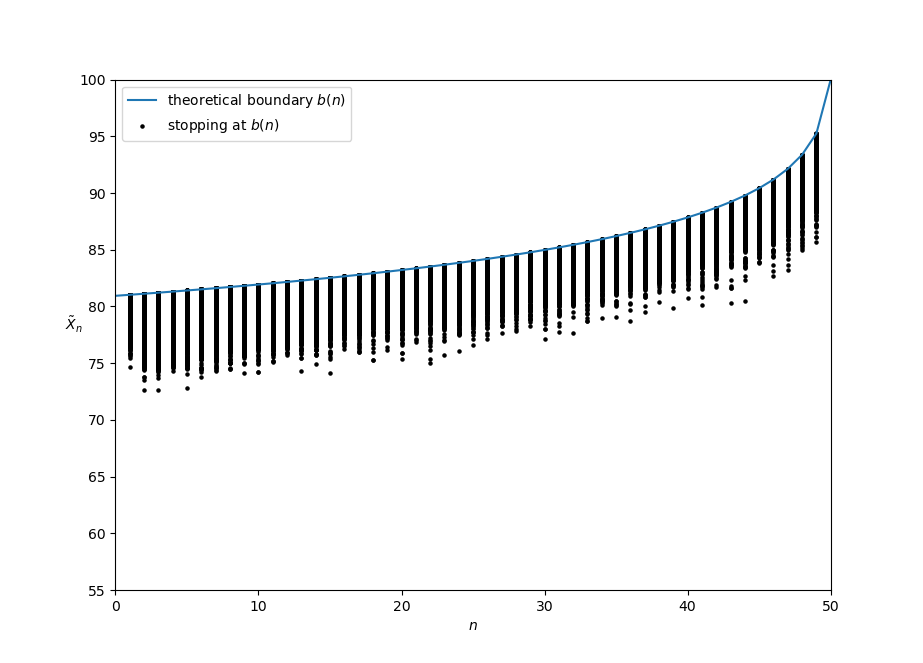}
\includegraphics[width=0.46\textwidth,height = 5cm]{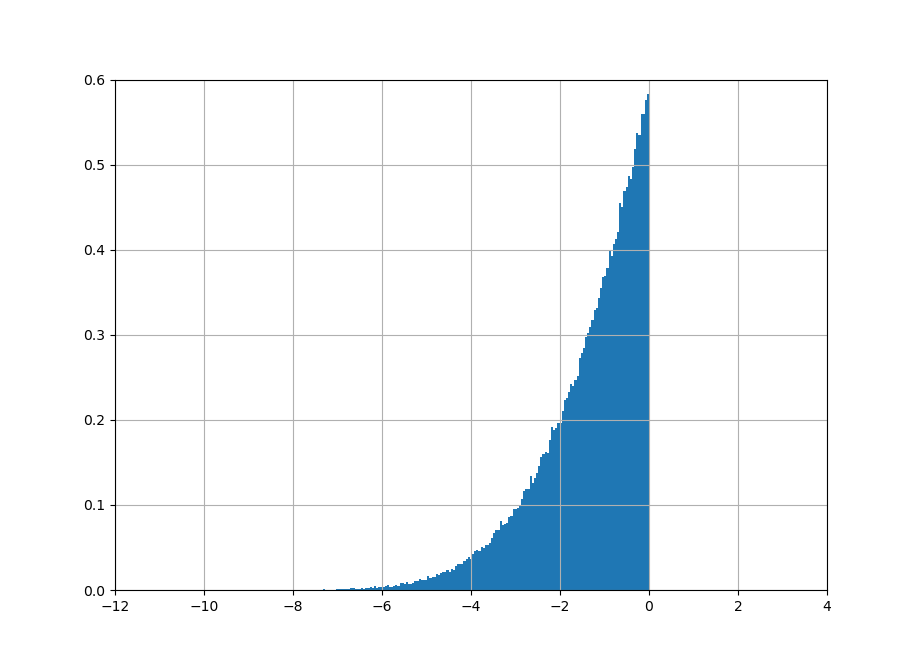}\\
{\bf\tiny Figure 7a}\hspace{5.5cm}{\bf\tiny Figure 8a}   
\end{center}

\begin{center}
\includegraphics[width=0.46\textwidth,height = 5cm]{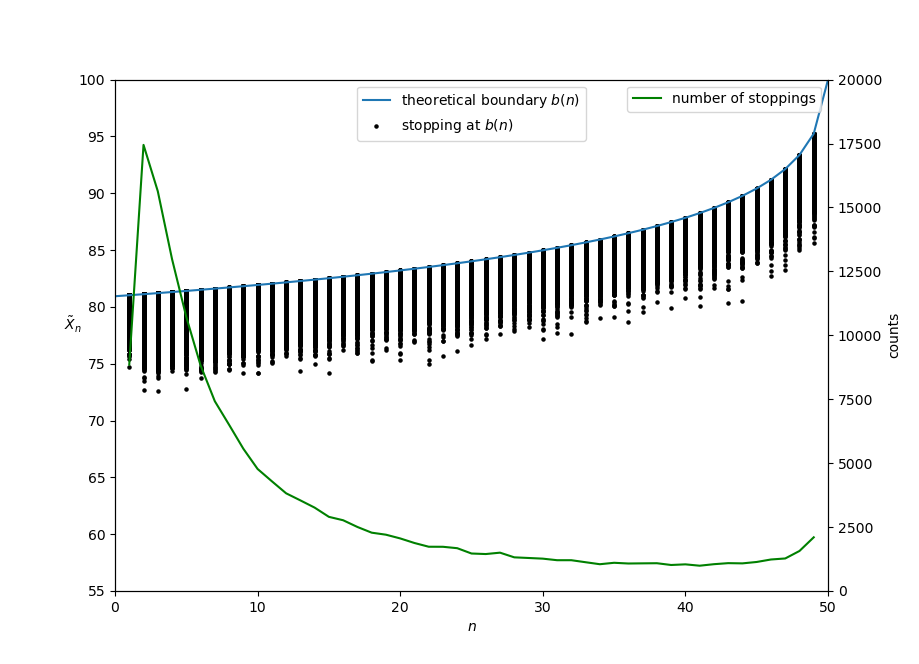}
\includegraphics[width=0.46\textwidth,height = 5cm]{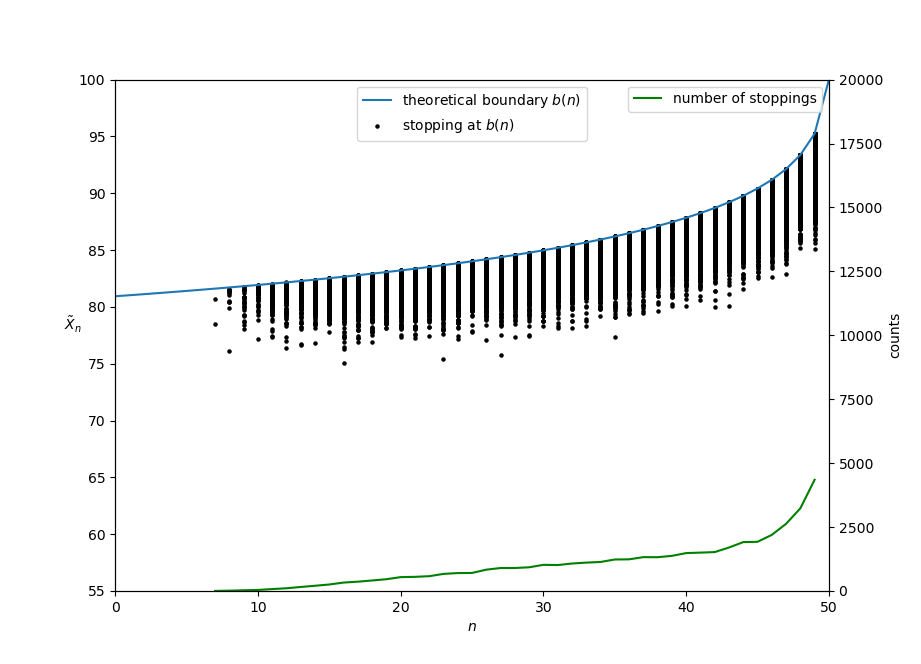}\\      
{\bf\tiny Figure 9a: ${\rm x}_M=85$}\hspace{4.5cm}{\bf\tiny Figure 10a: ${\rm x}_M=110$}
\end{center}

We substract $b(n)$ from any value associated with a bullet point
plotted underneath $b(n)$, which gives step-$n$-residuals, and Figure 8a shows the distribution
of all residuals, across all steps $n=1,\dots,N-1$.
These residuals are all negative because a continuous time stopping rule has been applied
to discretised geometric Brownian motion: 
if a test sample path $\tilde{\omega}^{(k)}$ is stopped at $1\le n\le N-1$
then $\tilde{\omega}^{(k)}_{n-1}>b(n)$ but $\tilde{\omega}^{(k)}_{n}<b(n)$.
As a consequence, 
$\tilde{V}_{0,0}$ is not attained at (\ref{notOptimal})---the stopping seems to happen too late.

In Figure 9a, the additional green chart gives 
the number of test sample paths which were stopped at $n=1,\dots,N-1$.
Since ${\rm x}_M=85$ is close to the theoretical boundary, many test sample paths
cross this boundary very soon. In contrast, in Figure 10a, where ${\rm x}_M=110$,
less test sample paths cross the theoretical boundary and they need longer to do so,
while many test sample paths would not cross it at all.

The bullet points in Figure 7b below refer to
stopping test sample paths with respect to the stopping rule obtained when plugging
$\tilde{X}_n,\,n=0,\dots,N$, into the right-hand side of (\ref{ourTauZero}),
where $\hat{g}_n,\,n=0,\dots,N$, were estimated using the training data for
our $\Delta$-algorithm combined with bagging and cross validation.
Again, the number of bags is $B=10$, 
and maximum tree-depth $(=10)$ as well as minimum node-size $(=10)$ were constrained, too.
\begin{center}
\includegraphics[width=0.46\textwidth,height = 5cm]{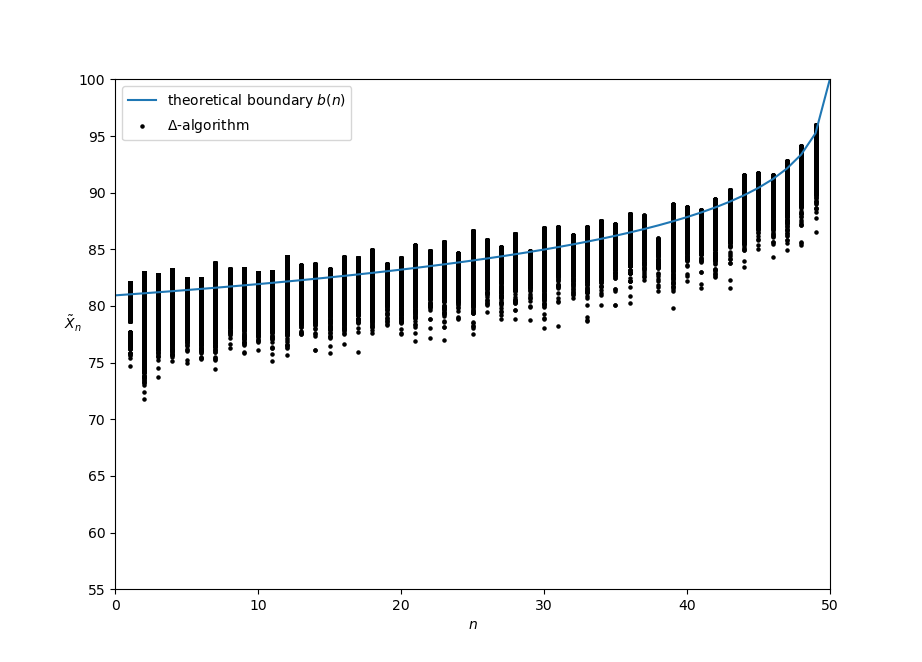}
\includegraphics[width=0.46\textwidth,height = 5cm]{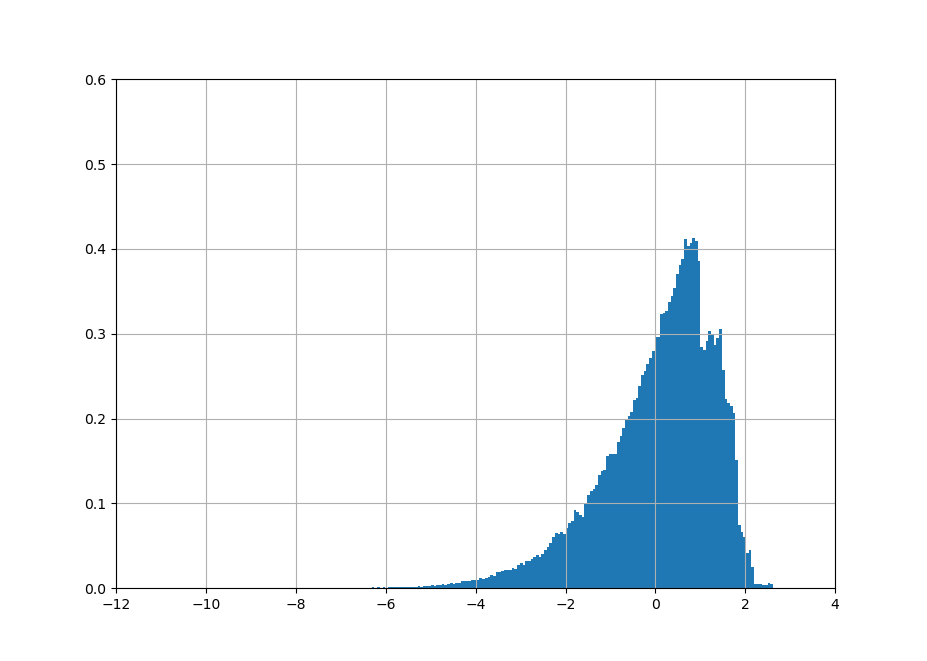}\\ 
{\bf\tiny Figure 7b}\hspace{5.5cm}{\bf\tiny Figure 8b}
\end{center}

\begin{center}\label{green charts}
\includegraphics[width=0.46\textwidth,height = 5cm]{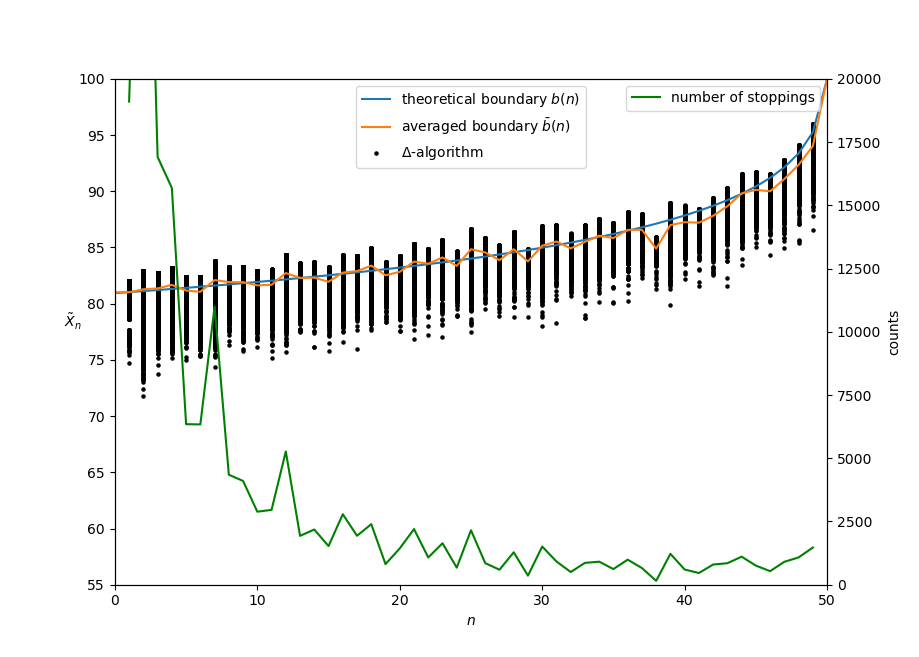}
\includegraphics[width=0.46\textwidth,height = 5cm]{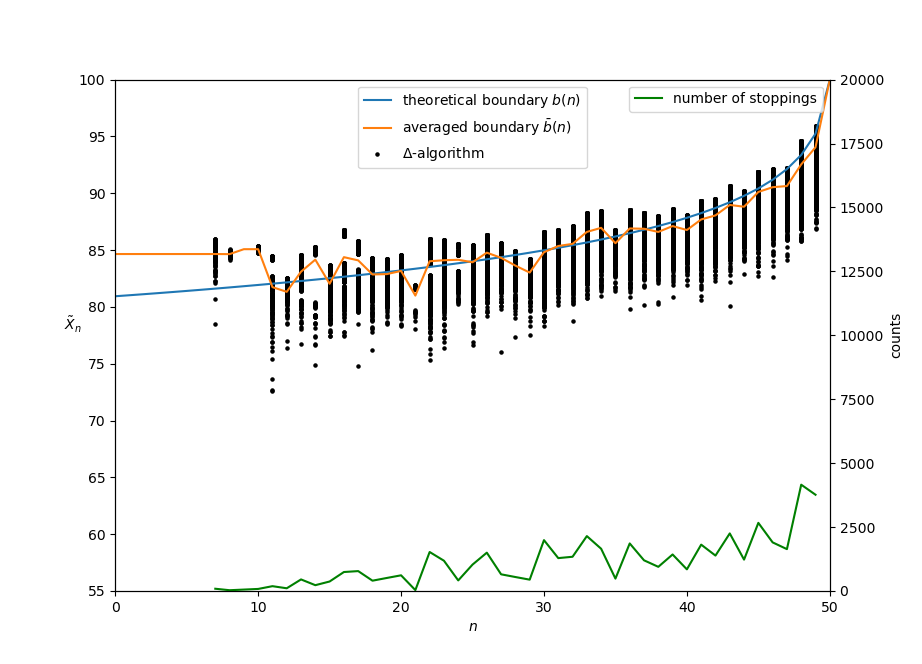}\\
{\bf\tiny Figure 9b: ${\rm x}_M=85$}\hspace{4.5cm}{\bf\tiny Figure 10b: ${\rm x}_M=110$}
\end{center}

Subject to being slightly less asymmetric,
the shape of the distribution of the residuals in Figure 8b roughly resembles the shape
of the distribution shown in Figure 8a, only the location seems to have moved to the right:
many residuals are positive, that is many test sample paths were stopped above
the theoretical boundary. However, in Figure 9b, the additional red chart gives
the average of the values associated with bullet points at each $n=1,\dots,N-1$,
and this averaged boundary $\bar{b}(n)$ nicely resembles the theoretical boundary
indicating good fit. Note that the spread of the bullet points per $n$ is natural
for sequential discrete-step problems, it is caused by the spread of the distribution
associated with the step-$n$-increments of the sequential data 
used for training and testing.

Furthermore, comparing the green charts in Figure 9a and Figure 9b reveals that
more test sample paths are stopped at earlier steps when stopping is based
on the $\Delta$-algorithm, which goes along well with our statement that stopping
according to (\ref{notOptimal}) would be too late to be optimal.
This difference is less pronounced when comparing Figures 10a and 10b,
where ${\rm x}_M=110$ and test sample paths are stopped later.
Also, the averaged boundary in Figure 10b shows a better fit to the theoretical boundary
only at later steps when more test sample paths are stopped, though many
test sample paths, actually about $75\%$, are NOT stopped before $N$ at all.
If ${\rm x}_M=85$ on the other hand, 
only about $12\%$ of all test sample paths are NOT stopped before $N$.
Nevertheless, the green chart in Figure 10b tells that, even in absolute numbers,
there are more test sample paths with ${\rm x}_M=110$ stopped at later steps
than test sample paths with ${\rm x}_M=85$. 
So, when using the $\Delta$-algorithm for plotting stopping boundaries, it should be beneficial 
to generate test sample paths with respect to a mixture of initial ${\bf x}_M$-values.

Choosing ${\rm x}_M=85$ in both cases, Figures 11 and 12 below shed some light on 
how the $\Delta$-algorithm behaves with respect to different sample sizes.
It can be clearly seen that the averaged boundary's fit to the theoretical boundary
is worse for $\tilde{K}=20000$. However, smoothening the red chart in Figure 12
should still give a curve very close to the theoretical boundary,
which might be sufficient for practitioners.
\begin{center}
\includegraphics[width=0.45\textwidth,height = 5cm]{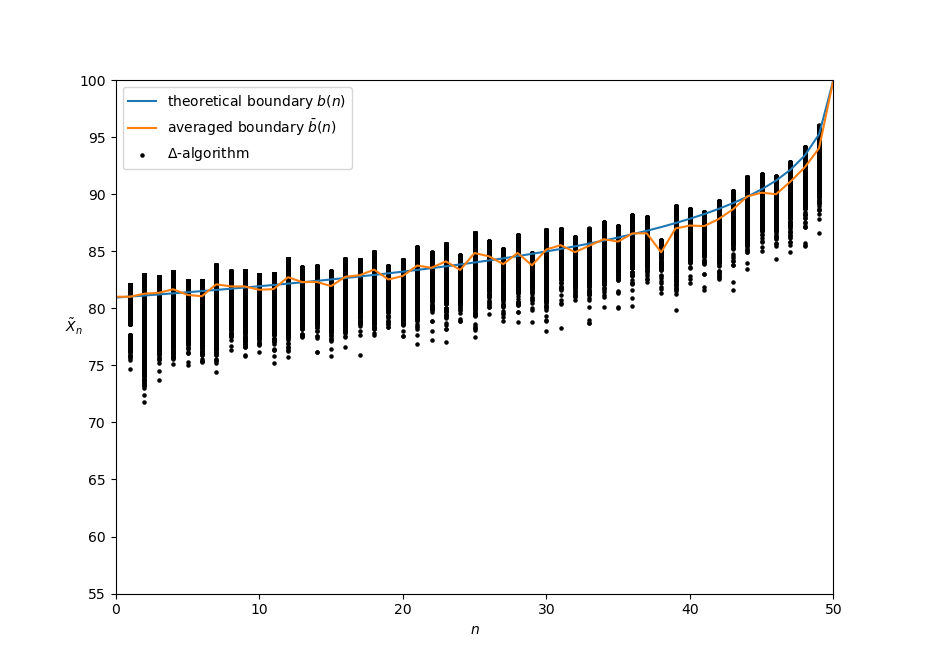} 
\hskip 5pt
\includegraphics[width=0.45\textwidth,height = 5cm]{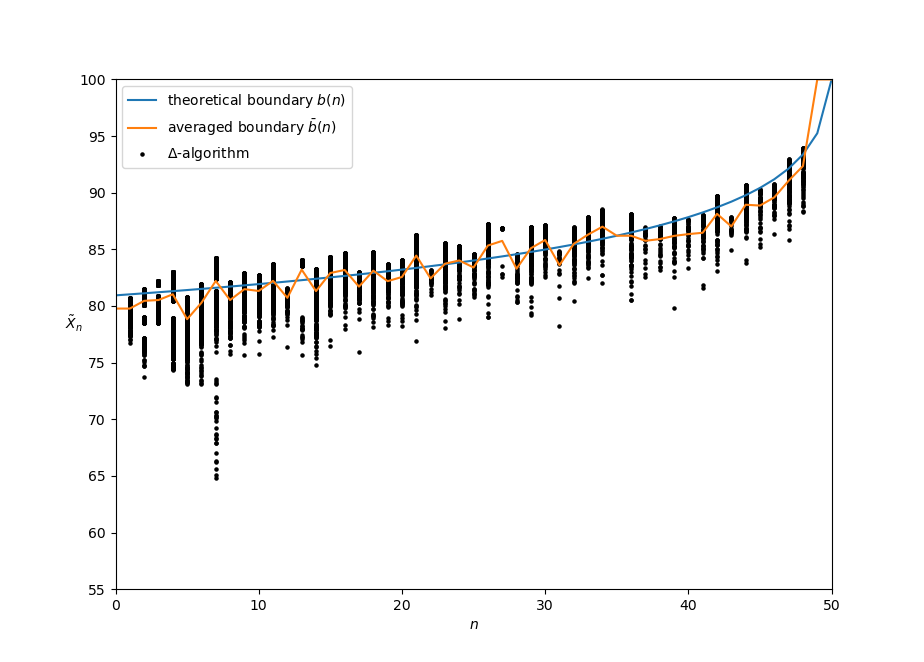}\\
{\bf\tiny Figure 11: $\tilde{K} = 200000$}\hspace{4.5cm}{\bf\tiny Figure 12: $\tilde{K} = 20000$}
\end{center}
\subsection{Pricing High-Dimensional Bermudan Max-Call Options}\label{maxCall1}
We are going to address two issues with respect to the results obtained in \cite[Section 4.1]{BCJ2019},
first, reproduction of these results,                                                                     second, comparing them with the results produced by our $\Delta$-algorithm.

The application treated in \cite[Section 4.1]{BCJ2019}
deals with the optimal stopping problem (\ref{trueValue})
associated with reward function
\begin{equation}\label{maxCallReward}
u(n,{\bf x})\,=\,
\exp\{-rnT/N\}(\max_{1\le d\le D}{\bf x}[d]-C\,)^+,
\quad n\in\{0,\dots,N\},\,{\bf x}\in\bfR^D,
\end{equation}
and random variables $X_n,\,n=0,\dots,N$, following a stochastic process in $\bfR^D$
whose components are independent discretised geometric Brownian motions, 
where the parameters $r,T,N$, and $C$, have got the same meaning as in Section \ref{theBoundary}.

As suggested in \cite{BCJ2019}, our training data has been generated via
\[
{\bf x}_{M+\,n}^{(k)}[d]
\,=\,
{\rm x}_M\times\exp\{\,
(\mu-\frac{\sigma[d]^2}{2})\frac{nT}{N}      
+\sigma[d]\sqrt{\frac{T}{N}}\sum_{n'=1}^n\varepsilon_{n'}^{(k)}[d]
\,\}
\]
using i.i.d.\ standard normals 
$\varepsilon_{n'}^{(k)}[d],\,n'=1,\dots,N,\,k=1,\dots,K,\,d=1,\dots,D$,
while our test data following the same dynamics
has been generated using another but independent batch of i.i.d.\ standard normals
$\tilde{\varepsilon}_{n'}^{(k)}[d],\,n'=1,\dots,N,\,k=1,\dots,\tilde{K},\,d=1,\dots,D$.
For $n=0,\dots,N$,
let $X_n$ and $\tilde{X}_n$ denote the corresponding canonical random variables.

For all numerical experiments in this subsection,
the interest rate $r=0.05$, drift $\mu=-0.05$, maturity $T=3$, number of time periods $N=9$,
and strike level $C=100$, will be fixed.
Otherwise, we consider the following two scenarios 
for different choices of initials ${\rm x_M}>0$, and dimensions $D$.\\

{\bf Symmetric Case:}
\[
\sigma[d]=0.2,\quad d=1,\dots,D.
\]

{\bf Asymmetric Case:}
\[
\sigma[d]\,=\left\{
\begin{array}{ccl}
0.08+0.32\times(d-1)/(D-1)&:&d=1,\dots,D\le 5,\\
0.1+d/(2D)&:&d=1,\dots,D>5.
\end{array}\right.
\]

Let a {\it feature} be any functional of the components of a vector ${\bf x}\in\bfR^D$,
which makes the components themselves features, too.
Then, for each $n=0,\dots,N$,
the step-$n$-reward $u(n,X_n)$ is a feature of the random variable $X_n$, 
and the $\Delta$-values crucial for training our CART-trees are entirely based on these features, only. 
This raises the question whether split-points could be calculated with respect to the values
of the step-$n$-rewards, only, making the problem artificially one-dimensional.
Unfortunately, in the majority of all cases, the sequence of random variables
$u(n,X_n),\,n=0,\dots,N$, would NOT form a Markov process with respect to the filtration
generated by this reward-process, and hence finding split-points this way would produce
rather inaccurate estimates of optimal stopping rules.

\label{features}
However, supplementing the reward-process with a few more features could create a still
quite low-dimensional stochastic process which is at least approximately\footnote{
The meaning of `approximation' in this context has to be clarified, of course.}
Markov enough for numerical purposes.
If this is the case then one could indeed choose split-components from the smaller
number of new features followed by calculating split-points with respect to these
features, only.

The problem with this approach is that, if they do exist,
it is NOT easy to find the few extra features
the reward-process has to be supplemented with,
at least we are not aware of any algorithmic solution. Simply experimenting,
we nevertheless found good extra features for the rewards given by (\ref{maxCallReward}),
and we will present results using these features further below.

Experimenting, though, can take a long time, thus, in practical terms,
split-components should be chosen from the components of the original Markov process
$X_n,\,n=0,\dots,N$. But, since working out $\Delta$-values requires knowledge
of the reward process, it would be good practice to add the step-$n$-rewards 
creating a $(D+1)$-dimensional Markov process, $(X_n,u(n,X_n)),\,n=0,\dots,N$, 
as then our CART-trees would not have to learn the structure of the reward function.

This `trick' has been mentioned and applied in \cite{BCJ2019}, and we will do the same.
We at first find the estimates $\hat{f}_{0,n},\,n=1,\dots,N$,
using $(X_n,\,u(n,X_n)),\,n=0,\dots,N$, for our $\Delta$-algorithm
combined with bagging and cross-validation, where the number of bags is $B=10$,
and maximum tree-depth (=10) as well as minimum node-size (=10) were constrained, too,
and we then calculate $V_{0,0}^{test}$ using the values of the test-random-variables
$(\tilde{X}_n,\,u(n,\tilde{X}_n)),\,n=0,\dots,N$,
where the number of sample paths generated for training and testing
has been chosen to be $K=100000$ and $\tilde{K}=4096000$, respectively.
The unusual large number of test-sample-paths is required to match
the number of test-sample-paths used in \cite[Section 4.1]{BCJ2019}.

To demonstrate the effectiveness of this `trick', we re-do the above experiment
only using $X_n,\,n=0,\dots,N$, and $\tilde{X}_n,\,n=0,\dots,N$,
for training and testing, respectively.

We observe that, 
except for $D=10$ and ${\rm x}_M=90$ in the Symmetric Case,
adding the step-$n$-rewards gives larger $V_{0,0}^{test}$-values in all dimensions $D\ge 10$ 
in both Tables \ref{sym with and without gain} and \ref{asym with and without gain}.
So, similar to the reasoning leading to Corollary \ref{better estimate},
adding step-$n$-rewards gives better $V_{0,0}^{test}$-values in higher dimensions.
Moreover, except for $D=500$ in both the Symmetric and Asymmetric Case,
the computing time is always shorter which we think has to do with our comment above that
the CART-trees would not have to learn the structure of the reward function.
At the moment we cannot explain why this principle would not apply in the
exceptional case of $D=500$.

Tables \ref{symComparing} and \ref{asymComparing} present the main results of this paper.
Here, both CART-trees and Neural Networks were trained on
$(D+1)$-dimensional values of training-random-variables, $({X}_n,u(n,{X}_n)),\,n=0,\dots,N$,
which is consistent with the approach taken in \cite{BCJ2019}.

Training CART-trees by applying our $\Delta$-algorithm 
combined with bagging ($B=10$) and cross-validation, 
while constraining maximum tree-depth (=10) as well as minimum node-size (=10),
yields estimates $\hat{f}_{0,n},\,n=1,\dots,N$,
which were used to calculate $V_{0,0}^{test}$ according to (\ref{ourLowerBound}),
plugging in test data of size $\tilde{K}=4096000$ as suggested in \cite{BCJ2019}.

Following \url{https://github.com/erraydin/Deep-Optimal-Stopping},
we implemented the algorithm suggested in \cite{BCJ2019} for training
Neural Networks which yield estimates ${f}^\star_{0,n},\,n=1,\dots,N$,
which were used to calculate $V_{0,0}^\star$ according to (\ref{theirLowerBound}),
plugging in the \underline{same} test data of size $\tilde{K}=4096000$ 
we had employed to find $V_{0,0}^{test}$.

Note that CART-trees and Neural Networks were trained on different training data.
For CART-trees, we generated training data of size $K=100000$,
while, to be consistent with \cite{BCJ2019}, depending on the dimension $D$,
the Neural Networks were trained on huge data of size $8192\times(3000+D)$,
which is about $240$ times the size of the CART-tree training data.

The computation times quoted in Tables \ref{symComparing} and \ref{asymComparing}
refer to computations on a machine with specifications:
Intel(R) Core(TM) i7-6700 CPU with 3.40GHz, 3408 Mhz,              
4 Core(s), 8 Logical Processor(s), and 16GB memory.
Since 16GB memory was too small for the huge Neural Network training data,
we split it generating fixed random seeds for each of the steps $n=1,\dots,N$,
which slowed down computation time for the corresponding estimates by about $2\%$.

Nevertheless, our computation times of $V^\star_{0,0}$ do not seem to match 
the computation times reported by the authors of \cite{BCJ2019} 
although the specifications of the machine they used do not differ that much from ours.
The reason might be that their calculations are more GPU-based\footnote{
By rule of thumb, GPU-calculations are about 30 times faster than CPU-calculations.}
than ours, or that their implementation is `smarter' in terms of exploiting hardware
than the rather basic implementation we followed.
Still, we would not have expected a factor as large as more than 100 times slower.

The values of $V^\star_{0,0}$, though, 
match the values reported by the authors of \cite{BCJ2019} extremely well,
and hence we can trust our implementation of their algorithm.
Below, still discussing Tables \ref{symComparing} and \ref{asymComparing},
we are going to compare $V^{test}_{0,0}$ and $V^\star_{0,0}$,
the latter of which is taken to be the point of reference.

To start with, we observe that, except for $D=500$ and $x_M=110$ in the Symmetric Case,
all values of $V^\star_{0,0}$ are bigger than those of $V^{test}_{0,0}$.
Thus, by Corollary \ref{better estimate}, 
the estimated Neural Networks seem to act better on unseen data than the estimated CART-trees.

However, the difference between $V^{test}_{0,0}$ and $V^\star_{0,0}$ is less pronounced 
in the Symmetric Case and becomes smaller and smaller the higher the dimension,
which can be clearly seen in the below figures,
where average relative errors have been plotted,
depending on $D$ and ${\rm x}_M$ as follows:
from all ${\rm x}_M$-related values of $|V^{test}_{0,0}-V^\star_{0,0}|/V^\star_{0,0}$\,,
take the average of the values associated with $D,\dots,500$.
\begin{center}  
\includegraphics[width=0.45\textwidth,height = 5cm]{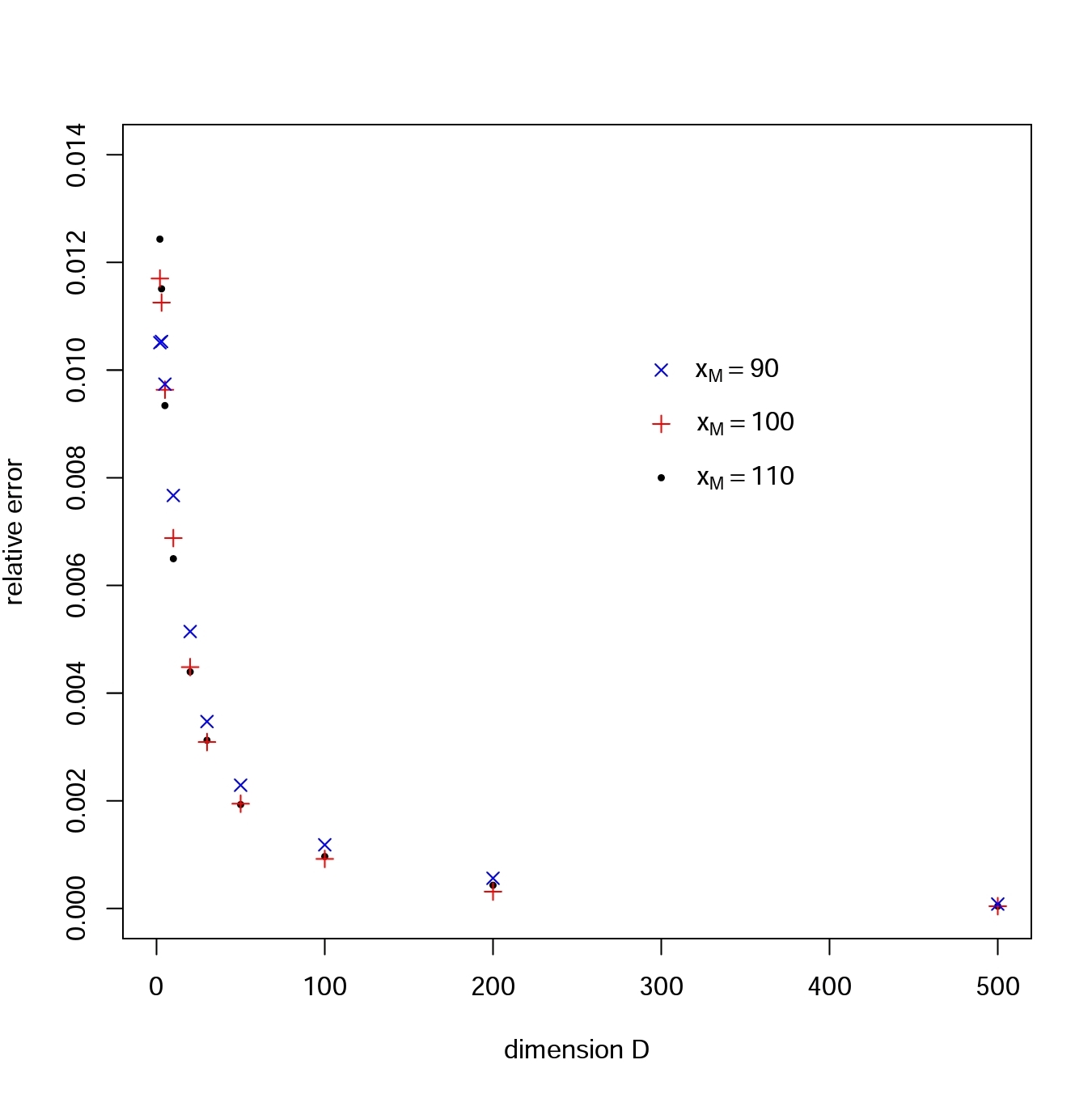}
\hskip 5pt             
\includegraphics[width=0.45\textwidth,height = 5cm]{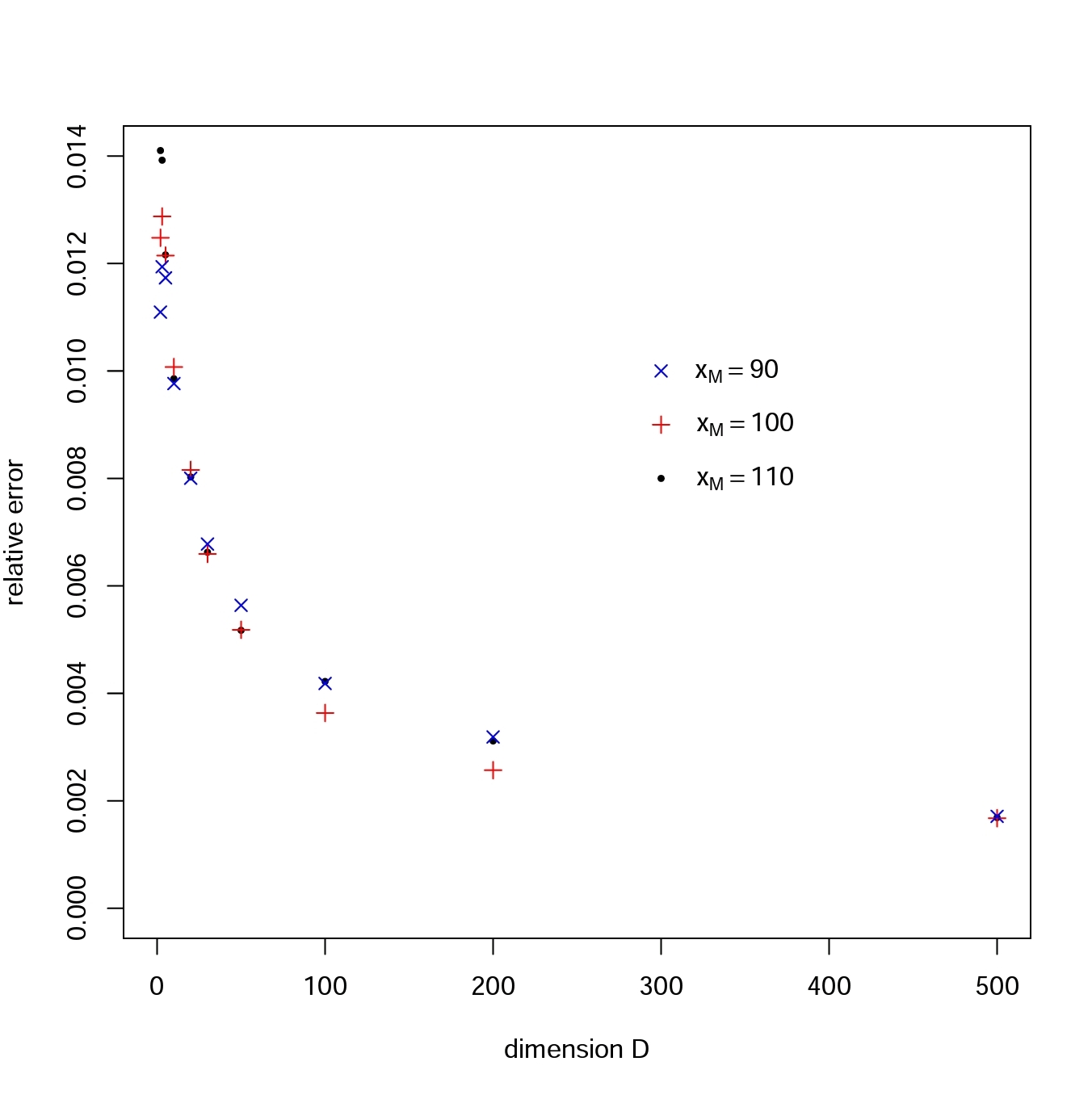}\\
$\rule{15pt}{0pt}$
{\bf\tiny Figure 13: $\mbox{\rm Symmetric Case}$}
\hspace{3.7cm}
{\bf\tiny Figure 14: $\mbox{\rm Asymmetric Case}$} 
\end{center}

Let us discuss the case $D=10$, where the averages have been
taken over relative errors associated with $D=10,20,30,50,100,200,500$.
Each figure shows three average relative errors at $D=10$,
which are spread between $0.6\%$ and $0.8\%$ in the Symmetric Case,
and which cluster at about $1\%$ in the Asymmetric Case.
But, why would at $D=10$ in the Symmetric Case
the average relative errors decrease with increasing values of ${\rm x}_M$?
\begin{remark}\rm\label{heuristic}
The key to answer the above question is the following obvious fact:
sample paths which only stop at maturity would only
contribute to the estimation of $\hat{f}_{0,N}$, while the estimation of 
the remaining $\hat{f}_{0,1},\dots,\hat{f}_{0,N-1}$ would be entirely based on
sample paths which stop before maturity. Thus, if less sample paths stop before maturity
then the number of sample paths available for estimating $\hat{f}_{0,1},\dots,\hat{f}_{0,N-1}$
might be low, which could leave some of these estimates with a larger variance.
\end{remark}

Now, recall the green charts in Figures 9b and 10b on page \pageref{green charts}.
There, when ${\rm x}_M$ was closer to the stopping boundary,
many more paths stopped before maturity.

But, ${\rm x}_M$ was one-dimensional in Figures 9b and 10b. 
In the present max-call example, we have got paths with initial values
\[
{\bf x}_M[d]\,=\,{\rm x}_M,\quad d=1,\dots,D.
\]

So, expecting the stopping behaviour to be similar to what we observed in Figures 9b and 10b,
the above remark suggests that 
the further away a multi-dimensional ${\bf x}_M$ is from the stopping boundary
the larger the variance of at least some of the estimates
$\hat{f}_{0,1},\dots,\hat{f}_{0,N-1}$,
which would eventually lead to a larger relative error of 
$V^{test}_{0,0}$ given by (\ref{ourLowerBound}).

Thus, if the initial values ${\bf x}_M$ associated with ${\rm x}_M=90,100,110$ 
are close enough to the stopping boundary in dimensions $D=10,20,30,50,100,200,500$,
and if ${\bf x}_M$ with larger ${\rm x}_M$ are closer to the stopping boundary in the Symmetric Case,
then the corresponding average relative errors should indeed be smaller,
as shown at $D=10$ in Figure 13.

We do not know how close to the stopping boundary 
the initial values ${\bf x}_M$ associated with ${\rm x}_M=90,100,110$ really are.
Nevertheless, in the Symmetric Case, 
we expect the stopping boundary to be symmetric in some sense, and we therefore think that
${\bf x}_M$ with larger ${\rm x}_M$ should be closer to the stopping boundary,
because, when ${\rm x}_M$ increases, the initial ${\bf x}_M$ move up on the diagonal in $\bfR^D$. 

Of course, the above approach is rather heuristic, but it would explain what we see
at $D=10$ in Figure 13. Moreover, it would also explain why
including relative errors associated with smaller dimensions, like $D=2,3,5$, 
would make the average relative errors less representative in the case of max-call options.
Indeed, if $D=1$, the max-call would be an ordinary American call, and, as such, 
its value would coincide with the value of the corresponding European call
with no early exercise opportunity, in other words, all paths would only stop at maturity.
Of course, if $D$ is too close to one, 
many paths would still stop at maturity, only, which, as explained in Remark \ref{heuristic}, 
would lead to larger variances of the estimates $\hat{f}_{0,1},\dots,\hat{f}_{0,N-1}$, 
and hence to larger relative errors of $V^{test}_{0,0}$.

The larger relative error of $V^{test}_{0,0}$ in smaller dimensions
is even more noticeable in Figure 14.
We also think that relative errors of $V^{test}_{0,0}$ are still larger
when the initial values ${\bf x}_M$ are further away from the stopping boundary.
In the Asymmetric Case, though,
increasing the values of ${\rm x}_M$, that is moving ${\bf x}_M$ up the diagonal in $\bfR^D$,
would not necessarily bring these initial values closer to the stopping boundary, 
which we expect to be less symmetric than in the Symmetric Case, by obvious reasons.
So, what we see in Figure 14 tells that, in the Asymmetric Case, 
it cannot be concluded from the order of ${\rm x}_M=90,100,110$
whether the corresponding initial values ${\bf x}_M$ are further away from or closer to
the stopping boundary.

Next, in Tables \ref{symComparing} and \ref{asymComparing},
we should not compare the computation times for $V^{test}_{0,0}$ and $V^\star_{0,0}$
as they include the times for training and the size of data used   
to find $f^\star_{0,n},\,n=0,\dots,N$, is about 240 times larger    
than the size of data we used to find $\hat{f}_{0,n},\,n=0,\dots,N$.
So, the big difference between the corresponding computation times is not surprising. 

The authors of \cite{BCJ2019} did not say why they needed training data of this size,
but it might be that they observed values of insufficient accuracy
when using smaller training data, as that is what we observed. 

For example, training the Neural Networks on data of size $K=100000$ produces 
the results shown in Tables \ref{symComparing_100000} and \ref{asymComparing_100000}.
At that size of training data, their algorithm does not seem to have converged, yet.
The values of $V^{test}_{0,0,}$ based on the $\Delta$-algorithm, using CART-trees 
trained on data of the same size, are much closer to the reference-values (last column),
suggesting that the $\Delta$-algorithm converges much faster.
Since the algorithm based on Neural Networks has not yet converged,
the much faster computation times reported in 
Tables \ref{symComparing_100000} and \ref{asymComparing_100000} are deceptive:
for Neural Networks to produce values as good as $V^{test}_{0,0}$, 
the size of the training data, $K$, has to be much larger than $100000$, leading to 
rather long computation times as reported in Tables \ref{symComparing} and \ref{asymComparing},
at least in dimensions $10\le D\le 100$. 

Interestingly, even in the case of smaller sized training data,
both algorithms perform exceptionally well in very high dimensions $D\ge 200$.
We therefore think that
the $\mbox{\rm max-call}$ example is not a hard example for estimating optimal stopping rules
in very high dimensions---it appears to be harder in medium-size dimensions $10\le D\le 100$.

Actually, when experimenting with features 
as explained on page \pageref{features} at the beginning of this section, 
we found that using the four-dimensional stochastic process
\[
(\,u(n,X_n),\,\phi_1(X_n),\,\phi_2(X_n),\,\phi_3(X_n)\,)^T,\quad n=0,\dots,N,
\]
for our $\Delta$-algorithm, where\footnote{Here, $d^\star\,=\,{\rm arg}\max_{1\le d\le D}x[d]$.}
\begin{equation}\label{ourPhi}
\phi_1({\bf x})\,=\,\mbox{$\max_{d}x[d]$},\quad
\phi_2({\bf x})\,=\,\mbox{$\max_{d\not=d^\star}x[d]$},\quad
\phi_3({\bf x})\,=\,\phi_1({\bf x})-\phi_2({\bf x}),
\end{equation}
produces excellent results---see
Tables \ref{symCompare4features} and \ref{asymCompare4features}. 
Again, the algorithm was  combined with bagging (number of bags $B=10$) and cross-validation,
also constraining tree-depth (=10) as well as minimum node-size (=10).

All values of $V_{0,0}^{test}$ are very close to the corresponding reference-values (last column).
But, when $D$ is large 
($D\ge 10$ in the Symmetric case and $D\ge 50$ in the Asymmetric Case),
$V_{0,0}^{test}$ even dominates $V^\star_{0,0}$,
and hence Corollary \ref{better estimate} implies that CART-trees 
trained on $100000$ sample paths of a four-dimensional process
act better on unseen data than Neural Networks trained on
$8192\times(3000+D)$ sample paths of a $(D+1)$-dimensional process.

So, the max-call example for estimating optimal stopping rules is not only easier
in very high dimensions, the information needed to effectively estimate the wanted
optimal stopping rules can be given by four features, only,
whatever the dimension. It would be interesting to know whether checking the weights
of the trained  Neural Networks would reveal anything about the four features 
we experimentally found.

The concluding remark of this section concerns the computation times reported in
Tables \ref{symCompare4features} and \ref{asymCompare4features}.
It is striking that training the CART-trees only takes about a tenth of the time
needed to calculate $V^{test}_{0,0}$, which can be explained as follows:
the size of the test data, $\tilde{K}=4096000$, 
is about $40$ times the size of the training data,
to be consistent with how the values of $V^\star_{0,0}$ were calculated.
While the training times are more or less stable, the times for calculating
$V^{test}_{0,0}$ eventually increase in very high dimensions, which is again due to
the difference between training and test data size combined with the effect of
only four features being evaluated. Interestingly, we also observe that training
and calculation times can be larger in low dimensions than in the
medium-size dimensions $10\le D\le 100$, which we think is caused by deeper CART-trees
being `grown' in low dimensions.

\begin{table}[ht]
\centering
\resizebox{\textwidth}{9cm}{
\begin{tabular}{cccccc}
\hline
{$D$} &
{${\rm x}_M$} &
{\begin{tabular}[c]{@{}c@{}}$X_n$ \&\ $u(n,X_n)$\\ \rule{0pt}{11pt} $V_{0,0}^{test}$ \end{tabular}} &
{\begin{tabular}[c]{@{}c@{}} computation\\ time \end{tabular}} &
{\begin{tabular}[c]{@{}c@{}}$X_n$ only\\ \rule{0pt}{11pt} $V_{0,0}^{test}$ \end{tabular}} &
{\begin{tabular}[c]{@{}c@{}}computation\\ time \end{tabular}} \\ \hline
2   & 90  & 7.971   & 424.7  & 7.993   & 471.6  \\
2   & 100 & 13.656  & 427.7  & 13.751  & 492.1  \\
2   & 110 & 20.877  & 423.8  & 20.969  & 510.8  \\
3   & 90  & 11.070  & 467.7  & 11.111  & 592.7  \\
3   & 100 & 18.220  & 488.6  & 18.258  & 673.9  \\
3   & 110 & 26.755  & 551.3  & 26.856  & 696.2  \\
5   & 90  & 16.203  & 591.2  & 16.208  & 859.7  \\
5   & 100 & 25.350  & 622.3  & 25.408  & 906.6  \\
5   & 110 & 35.648  & 689.0  & 35.681  & 954.5  \\
10  & 90  & 25.554  & 833.7  & 25.444  & 1210.1 \\
10  & 100 & 37.389  & 827.4  & 37.268  & 1240.6 \\
10  & 110 & 49.753  & 864.8  & 49.559  & 1277.4 \\
20  & 90  & 37.117  & 1030.4 & 36.420  & 1434.8 \\
20  & 100 & 50.890  & 1030.8 & 49.990  & 1432.2 \\
20  & 110 & 64.709  & 1115.0 & 63.703  & 1489.0 \\
30  & 90  & 44.356  & 1053.8 & 43.255  & 1572.5 \\
30  & 100 & 58.952  & 1108.8 & 57.729  & 1596.5 \\
30  & 110 & 73.548  & 1174.4 & 72.071  & 1641.8 \\
50  & 90  & 53.555  & 1362.3 & 52.620  & 1735.9 \\
50  & 100 & 69.201  & 1369.9 & 68.026  & 1780.0 \\
50  & 110 & 84.800  & 1282.0 & 83.442  & 1824.6 \\
100 & 90  & 66.136  & 1426.5 & 65.327  & 2002.3 \\
100 & 100 & 83.139  & 1443.7 & 82.148  & 2069.1 \\
100 & 110 & 100.119 & 1579.6 & 98.969  & 2184.2 \\
200 & 90  & 78.824  & 1998.1 & 78.127  & 2481.3 \\
200 & 100 & 97.236  & 1978.3 & 96.371  & 2447.8 \\
200 & 110 & 115.639 & 2044.4 & 114.615 & 2736.6 \\
500 & 90  & 95.869  & 3308.7 & 95.253  & 3194.5 \\
500 & 100 & 116.141 & 3546.5 & 115.400 & 2999.1 \\
500 & 110 & 136.426 & 3394.9 & 135.547 & 3775.6 \\  \hline
\end{tabular}
}
\caption{Results for the Symmetric Case, where computation time refers to the time
in seconds for finding the estimates $\hat{f}_{0,n},\,n=1,\dots,N$,
and calculating $V_{0,0}^{test}$.}
\label{sym with and without gain}
\end{table}

\begin{table}[ht]
\centering
\resizebox{\textwidth}{9cm}{
\begin{tabular}{cccccc}
\hline
$D$ &
{${\rm x}_M$} &
{\begin{tabular}[c]{@{}c@{}}$X_n$ \&\ $u(n,X_n)$\\ \rule{0pt}{11pt} $V_{0,0}^{test}$ \end{tabular}} &
{\begin{tabular}[c]{@{}c@{}} computation\\ time \end{tabular}} &
{\begin{tabular}[c]{@{}c@{}}$X_n$ only\\ \rule{0pt}{11pt} $V_{0,0}^{test}$ \end{tabular}} &
{\begin{tabular}[c]{@{}c@{}}computation\\ time \end{tabular}} \\ \hline
2   & 90  & 14.250  & 376.4  & 14.249  & 430.6  \\
2   & 100 & 19.581  & 440.0  & 19.570  & 510.7  \\
2   & 110 & 26.679  & 496.5  & 26.697  & 568.6  \\
3   & 90  & 18.793  & 495.8  & 18.804  & 591.7  \\
3   & 100 & 26.144  & 531.8  & 26.100  & 623.8  \\
3   & 110 & 34.799  & 592.6  & 34.678  & 701.8  \\
5   & 90  & 26.905  & 616.1  & 26.884  & 807.0  \\
5   & 100 & 36.928  & 640.6  & 36.934  & 815.3  \\
5   & 110 & 48.017  & 662.8  & 47.897  & 882.7  \\
10  & 90  & 84.040  & 737.4  & 83.629  & 990.8  \\
10  & 100 & 102.261 & 729.0  & 101.891 & 960.5  \\
10  & 110 & 120.962 & 779.2  & 120.466 & 1031.3 \\
20  & 90  & 124.065 & 964.7  & 123.247 & 1213.1 \\
20  & 100 & 147.089 & 931.3  & 146.430 & 1124.1 \\
20  & 110 & 170.594 & 1020.5 & 169.598 & 1222.2 \\
30  & 90  & 152.650 & 1054.9 & 151.680 & 1397.2 \\
30  & 100 & 179.132 & 1174.7 & 178.130 & 1445.0 \\
30  & 110 & 205.629 & 1110.1 & 204.414 & 1339.6 \\
50  & 90  & 194.132 & 1327.6 & 192.847 & 1717.9 \\
50  & 100 & 225.321 & 1294.6 & 223.934 & 1589.5 \\
50  & 110 & 256.832 & 1324.7 & 255.048 & 1656.7 \\
100 & 90  & 261.600 & 1958.3 & 259.958 & 2290.7 \\
100 & 100 & 300.279 & 1754.5 & 298.497 & 2119.2 \\
100 & 110 & 338.705 & 1883.2 & 336.870 & 2223.8 \\
200 & 90  & 342.817 & 2567.0 & 341.324 & 2974.9 \\
200 & 100 & 390.694 & 2470.5 & 388.815 & 2856.2 \\
200 & 110 & 438.116 & 2488.9 & 436.314 & 2953.8 \\
500 & 90  & 474.916 & 4241.8 & 473.289 & 4425.9 \\
500 & 100 & 537.467 & 4286.3 & 535.441 & 4153.6 \\
500 & 110 & 599.585 & 4645.3 & 597.590 & 4537.4 \\ \hline
\end{tabular}
}
\caption{Results for the Asymmetric Case, where computation time refers to the time
in seconds for finding the estimates $\hat{f}_{0,n},\,n=1,\dots,N$,
and calculating $V_{0,0}^{test}$.}
\label{asym with and without gain}
\end{table}

\begin{table}[ht]
\centering
\resizebox{\textwidth}{9cm}{
\begin{tabular}{ccccccc}
\hline
{$D$} &
{${\rm x}_M$} &
{\begin{tabular}[c]{@{}c@{}}CART-trees \\ \rule{0pt}{11pt} $V_{0,0}^{test}$ \\
\rule{0pt}{11pt} $K_{\rm train}=100000$  \end{tabular}} &
{\begin{tabular}[c]{@{}c@{}} computation \\ time \end{tabular}} &
{\begin{tabular}[c]{@{}c@{}} Neural Networks \\ \rule{0pt}{11pt} $V_{0,0}^\star$ \\
\rule{0pt}{11pt} $K_{\rm train}>24\times 10^6$ \end{tabular}} &
{\begin{tabular}[c]{@{}c@{}} computation \\ time \end{tabular}} &
{\begin{tabular}[c]{@{}c@{}} reported \\ in \cite{BCJ2019} \end{tabular}} \\ \hline
2   & 90  & 7.971   & 424.7  & 8.054   & 924.0   & 8.072 (28.7)    \\
2   & 100 & 13.656  & 427.7  & 13.874  & 922.4   & 13.895 (28.7)   \\
2   & 110 & 20.877  & 423.8  & 21.319  & 913.1   & 21.353 (28.4)   \\
3   & 90  & 11.070  & 467.7  & 11.260  & 1035.4  & 11.290 (28.8)   \\
3   & 100 & 18.220  & 488.6  & 18.672  & 1018.6  & 18.690 (28.9)   \\
3   & 110 & 26.755  & 551.3  & 27.550  & 1016.5  & 27.564 (27.6)   \\
5   & 90  & 16.203  & 591.2  & 16.605  & 1258.9  & 16.648 (27.6)   \\
5   & 100 & 25.350  & 622.3  & 26.105  & 1267.7  & 26.156 (28.1)   \\
5   & 110 & 35.648  & 689.0  & 36.722  & 1252.9  & 36.766 (27.7)   \\
10  & 90  & 25.554  & 833.7  & 26.151  & 1779.0  & 26.208 (30.4)   \\
10  & 100 & 37.389  & 827.4  & 38.201  & 1780.9  & 38.321 (30.5)   \\
10  & 110 & 49.753  & 864.8  & 50.722  & 1777.9  & 50.857 (30.8)   \\
20  & 90  & 37.117  & 1030.4 & 37.625  & 2751.2  & 37.701 (37.2)   \\
20  & 100 & 50.890  & 1030.8 & 51.479  & 2751.4  & 51.571 (37.5)   \\
20  & 110 & 64.709  & 1115.0 & 65.412  & 2783.8  & 65.494 (37.3)   \\
30  & 90  & 44.356  & 1053.8 & 44.723  & 4489.2  & 44.797 (45.1)   \\
30  & 100 & 58.952  & 1108.8 & 59.408  & 4887.6  & 59.498 (45.5)   \\
30  & 110 & 73.548  & 1174.4 & 74.134  & 4257.3  & 74.221 (45.3)   \\
50  & 90  & 53.555  & 1362.3 & 53.857  & 5863.7  & 53.903 (58.7)   \\
50  & 100 & 69.201  & 1369.9 & 69.550  & 5911.9  & 69.582 (59.1)   \\
50  & 110 & 84.800  & 1282.0 & 85.211  & 5938.1  & 85.229 (59.0)   \\
100 & 90  & 66.136  & 1426.5 & 66.297  & 10790.1 & 66.342 (95.5)   \\
100 & 100 & 83.139  & 1443.7 & 83.317  & 10753.5 & 83.380 (95.9)   \\
100 & 110 & 100.119 & 1579.6 & 100.323 & 10725.0 & 100.420 (95.4)  \\
200 & 90  & 78.824  & 1998.1 & 78.906  & 21488.3 & 78.993 (170.9)  \\
200 & 100 & 97.236  & 1978.3 & 97.293  & 21489.3 & 97.405 (170.1)  \\
200 & 110 & 115.639 & 2044.4 & 115.734 & 21492.6 & 115.800 (170.6) \\
500 & 90  & 95.869  & 3308.7 & 95.877  & 59728.3 & 95.956 (493.4)  \\
500 & 100 & 116.141 & 3546.5 & 116.146 & 59360.7 & 116.235 (493.5) \\
500 & 110 & 136.426 & 3394.9 & 136.420 & 59629.4 & 136.547 (493.7) \\  \hline
\end{tabular}
}
\caption{
Comparing the Symmetric Case: the first computation time refers to the time
in seconds for finding the estimates $\hat{f}_{0,n},\,n=1,\dots,N$,
and calculating $V_{0,0}^{test}$, when using the $\Delta$-algorithm,
while the second computation time refers to the time
in seconds for finding the estimates ${f}^\star_{0,n},\,n=1,\dots,N$,
and calculating $V_{0,0}^\star$, when using Neural Networks.
The last column shows the lower bound and the corresponding computation time (in brackets)                as reported in \cite[page 16]{BCJ2019}.
}
\label{symComparing}
\end{table}

\begin{table}[ht]
\centering
\resizebox{\textwidth}{9cm}{
\begin{tabular}{ccccccc}
\hline
{$D$} &
{${\rm x}_M$} &
{\begin{tabular}[c]{@{}c@{}}CART-trees \\ \rule{0pt}{11pt} $V_{0,0}^{test}$ \\
\rule{0pt}{11pt} $K_{\rm train}=100000$  \end{tabular}} &
{\begin{tabular}[c]{@{}c@{}} computation \\ time \end{tabular}} &
{\begin{tabular}[c]{@{}c@{}} Neural Networks \\ \rule{0pt}{11pt} $V_{0,0}^\star$ \\
\rule{0pt}{11pt} $K_{\rm train}>24\times 10^6$  \end{tabular}} &
{\begin{tabular}[c]{@{}c@{}} computation \\ time \end{tabular}} &
{\begin{tabular}[c]{@{}c@{}} reported \\ in \cite{BCJ2019} \end{tabular}} \\ \hline
2   & 90  & 14.250  & 376.4  & 14.300  & 1026.6  & 14.325 (26.8)   \\
2   & 100 & 19.581  & 440.0  & 19.757  & 1030.9  & 19.802 (27.0)   \\
2   & 110 & 26.679  & 496.5  & 27.105  & 1025.6  & 27.170 (26.5)   \\
3   & 90  & 18.793  & 495.8  & 19.052  & 1168.5  & 19.093 (26.8)   \\
3   & 100 & 26.144  & 531.8  & 26.642  & 1160.4  & 26.680 (27.5)   \\
3   & 110 & 34.799  & 592.6  & 35.802  & 1161.1  & 35.842 (26.5)   \\
5   & 90  & 26.905  & 616.1  & 27.609  & 1430.2  & 27.662 (28.0)   \\
5   & 100 & 36.928  & 640.6  & 37.940  & 1421.7  & 37.976 (27.5)   \\
5   & 110 & 48.017  & 662.8  & 49.415  & 1420.3  & 49.485 (28.2)   \\
10  & 90  & 84.040  & 737.4  & 85.785  & 2057.9  & 85.937 (31.8)   \\
10  & 100 & 102.261 & 729.0  & 104.514 & 2081.7  & 104.692 (30.9)  \\
10  & 110 & 120.962 & 779.2  & 123.535 & 2063.2  & 123.668 (31.0)  \\
20  & 90  & 124.065 & 964.7  & 125.842 & 3266.5  & 125.916 (38.4)  \\
20  & 100 & 147.089 & 931.3  & 149.477 & 3278.2  & 149.587 (38.2)  \\
20  & 110 & 170.594 & 1020.5 & 173.197 & 3288.3  & 173.262 (38.4)  \\
30  & 90  & 152.650 & 1054.9 & 154.401 & 4567.0  & 154.486 (46.5)  \\
30  & 100 & 179.132 & 1174.7 & 181.352 & 4550.6  & 181.275 (46.4)  \\
30  & 110 & 205.629 & 1110.1 & 208.219 & 4573.3  & 208.223 (46.4)  \\
50  & 90  & 194.132 & 1327.6 & 196.093 & 6953.7  & 195.918 (60.7)  \\
50  & 100 & 225.321 & 1294.6 & 227.557 & 6859.5  & 227.386 (60.7)  \\
50  & 110 & 256.832 & 1324.7 & 258.910 & 6876.5  & 258.813 (60.7)  \\
100 & 90  & 261.600 & 1958.3 & 263.226 & 12439.8 & 263.193 (98.5)  \\
100 & 100 & 300.279 & 1754.5 & 302.020 & 12585.4 & 302.090 (98.2)  \\
100 & 110 & 338.705 & 1883.2 & 340.899 & 12609.2 & 340.763 (97.8)  \\
200 & 90  & 342.817 & 2567.0 & 344.424 & 24860.6 & 344.575 (175.4) \\
200 & 100 & 390.694 & 2470.5 & 392.052 & 24986.6 & 392.193 (175.1) \\
200 & 110 & 438.116 & 2488.9 & 440.115 & 24803.5 & 440.037 (175.1) \\
500 & 90  & 474.916 & 4241.8 & 475.735 & 67526.1 & 476.293 (504.5) \\
500 & 100 & 537.467 & 4286.3 & 538.375 & 67346.3 & 538.748 (504.6) \\
500 & 110 & 599.585 & 4645.3 & 600.599 & 67646.9 & 601.261 (504.9) \\ \hline
\end{tabular}
}
\caption{
Comparing the Asymmetric Case: the first computation time refers to the time
in seconds for finding the estimates $\hat{f}_{0,n},\,n=1,\dots,N$,
and calculating $V_{0,0}^{test}$, when using the $\Delta$-algorithm,
while the second computation time refers to the time
in seconds for finding the estimates ${f}^\star_{0,n},\,n=1,\dots,N$,
and calculating $V_{0,0}^\star$, when using Neural Networks.
The last column shows the lower bound and the corresponding computation time (in brackets)
as reported in \cite[page 17]{BCJ2019}.
}
\label{asymComparing}
\end{table}

\begin{table}[ht]
\centering
\resizebox{\textwidth}{9cm}{
\begin{tabular}{ccccccc}
\hline
{$D$} &
{${\rm x}_M$} &
{\begin{tabular}[c]{@{}c@{}}CART-trees \\ \rule{0pt}{11pt} $V_{0,0}^{test}$ \\
\rule{0pt}{11pt} $K_{\rm train}=100000$  \end{tabular}} &
{\begin{tabular}[c]{@{}c@{}} computation \\ time \end{tabular}} &
{\begin{tabular}[c]{@{}c@{}} Neural Networks \\ \rule{0pt}{11pt} $V_{0,0}^\star$ \\ 
\rule{0pt}{11pt} $K_{\rm train}=100000$   \end{tabular}} &
{\begin{tabular}[c]{@{}c@{}} computation \\ time \end{tabular}} &
{\begin{tabular}[c]{@{}c@{}} Neurale Networks \\ $V^\star_{0,0}$ \\
\rule{0pt}{11pt} $K_{\rm train}>24\times 10^6$  \end{tabular}} \\ \hline
2   & 90  & 7.971   & 424.7  & 7.534   & 20.7  & 8.054 \\
2   & 100 & 13.656  & 427.7  & 9.755   & 15.1  & 13.874 \\
2   & 110 & 20.877  & 423.8  & 18.516  & 15.4  & 21.319 \\
3   & 90  & 11.070  & 467.7  & 9.995   & 17.9  & 11.260 \\
3   & 100 & 18.220  & 488.6  & 13.720  & 16.3  & 18.672 \\
3   & 110 & 26.755  & 551.3  & 19.919  & 15.9  & 27.550 \\
5   & 90  & 16.203  & 591.2  & 10.151  & 18.0  & 16.605 \\
5   & 100 & 25.350  & 622.3  & 15.698  & 18.7  & 26.105 \\
5   & 110 & 35.648  & 689.0  & 26.120  & 18.8  & 26.722 \\
10  & 90  & 25.554  & 833.7  & 10.768  & 25.5  & 26.151 \\
10  & 100 & 37.389  & 827.4  & 23.961  & 25.0  & 38.201 \\
10  & 110 & 49.753  & 864.8  & 28.206  & 23.9  & 50.722 \\
20  & 90  & 37.117  & 1030.4 & 9.552   & 41.5  & 37.625 \\
20  & 100 & 50.890  & 1030.8 & 29.271  & 41.8  & 51.479 \\
20  & 110 & 64.709  & 1115.0 & 33.061  & 41.8  & 65.412 \\
30  & 90  & 44.356  & 1053.8 & 12.083  & 55.7  & 44.723 \\
30  & 100 & 58.952  & 1108.8 & 31.301  & 54.5  & 59.408 \\
30  & 110 & 73.548  & 1174.4 & 49.760  & 54.7  & 74.134 \\
50  & 90  & 53.555  & 1362.3 & 27.634  & 76.3  & 53.857 \\
50  & 100 & 69.201  & 1369.9 & 26.850  & 75.1  & 69.550 \\
50  & 110 & 84.800  & 1282.0 & 42.775  & 76.5  & 85.211 \\
100 & 90  & 66.136  & 1426.5 & 65.696  & 132.1 & 66.297 \\
100 & 100 & 83.139  & 1443.7 & 82.279  & 132.9 & 83.317 \\
100 & 110 & 100.119 & 1579.6 & 76.893  & 132.6 & 100.323 \\
200 & 90  & 78.824  & 1998.1 & 78.352  & 239.9 & 78.906 \\
200 & 100 & 97.236  & 1978.3 & 96.489  & 242.2 & 97.293 \\
200 & 110 & 115.639 & 2044.4 & 114.840 & 239.9 & 115.734 \\
500 & 90  & 95.869  & 3308.7 & 95.258  & 565.5 & 95.877 \\
500 & 100 & 116.141 & 3546.5 & 115.407 & 589.4 & 116.146 \\
500 & 110 & 136.426 & 3394.9 & 135.555 & 588.6 & 136.420 \\
\hline
\end{tabular}
}
\caption{
Comparing the Symmetric Case: the first computation time refers to the time
in seconds for finding the estimates $\hat{f}_{0,n},\,n=1,\dots,N$,
and calculating $V_{0,0}^{test}$, when using the $\Delta$-algorithm,
while the second computation time refers to the time
in seconds for finding the estimates ${f}^\star_{0,n},\,n=1,\dots,N$, 
and calculating $V_{0,0}^\star$, when using Neural Networks with $K_{\rm train}=100000$.
}
\label{symComparing_100000}
\end{table}

\begin{table}[ht]
\centering
\resizebox{\textwidth}{9cm}{
\begin{tabular}{ccccccc}
\hline
{$D$} &
{${\rm x}_M$} &
{\begin{tabular}[c]{@{}c@{}}CART-trees \\ \rule{0pt}{11pt} $V_{0,0}^{test}$ \\ 
\rule{0pt}{11pt} $K_{\rm train}=100000$  \end{tabular}} &
{\begin{tabular}[c]{@{}c@{}} computation \\ time \end{tabular}} &
{\begin{tabular}[c]{@{}c@{}} Neural Networks \\ \rule{0pt}{11pt} $V_{0,0}^\star$ \\ 
\rule{0pt}{11pt} $K_{\rm train}=100000$   \end{tabular}} &
{\begin{tabular}[c]{@{}c@{}} computation \\ time \end{tabular}} &
{\begin{tabular}[c]{@{}c@{}} Neural Networks \\ $V^\star_{0,0}$ \\ 
\rule{0pt}{11pt} $K_{\rm train}>24\times 10^6$   \end{tabular}} \\ \hline
2   & 90  & 14.250  & 376.4  & 12.561  & 18.2  & 14.300 \\
2   & 100 & 19.581  & 440.0  & 13.366  & 22.7  & 19.757 \\
2   & 110 & 26.679  & 496.5  & 25.832  & 28.5  & 27.105 \\
3   & 90  & 18.793  & 495.8  & 16.334  & 29.5  & 19.052 \\
3   & 100 & 26.144  & 531.8  & 23.120  & 20.3  & 26.642 \\
3   & 110 & 34.799  & 592.6  & 28.655  & 18.6  & 35.802 \\
5   & 90  & 26.905  & 616.1  & 19.926  & 20.6  & 27.609 \\
5   & 100 & 36.928  & 640.6  & 27.938  & 19.1  & 37.940 \\
5   & 110 & 48.017  & 662.8  & 40.072  & 19.1  & 49.415 \\
10  & 90  & 84.040  & 737.4  & 33.894  & 24.5  & 85.785 \\
10  & 100 & 102.261 & 729.0  & 79.245  & 25.3  & 104.514 \\
10  & 110 & 120.962 & 779.2  & 66.470  & 24.8  & 123.535 \\
20  & 90  & 124.065 & 964.7  & 60.874  & 40.2  & 125.842 \\
20  & 100 & 147.089 & 931.3  & 72.936  & 44.9  & 149.477 \\
20  & 110 & 170.594 & 1020.5 & 96.568  & 44.3  & 173.197 \\
30  & 90  & 152.650 & 1054.9 & 63.810  & 61.1  & 154.401 \\
30  & 100 & 179.132 & 1174.7 & 68.447  & 60.3  & 181.352 \\
30  & 110 & 205.629 & 1110.1 & 89.107  & 58.6  & 208.219 \\
50  & 90  & 194.132 & 1327.6 & 50.282  & 88.5  & 196.093 \\
50  & 100 & 225.321 & 1294.6 & 166.935 & 91.1  & 227.557 \\
50  & 110 & 256.832 & 1324.7 & 112.522 & 91.7  & 258.910 \\
100 & 90  & 261.600 & 1958.3 & 247.726 & 154.4 & 263.226 \\
100 & 100 & 300.279 & 1754.5 & 289.766 & 156.2 & 302.020 \\
100 & 110 & 338.705 & 1883.2 & 323.216 & 159.4 & 340.899 \\
200 & 90  & 342.817 & 2567.0 & 342.824 & 275.3 & 344.424 \\
200 & 100 & 390.694 & 2470.5 & 390.136 & 281.3 & 392.052 \\
200 & 110 & 438.116 & 2488.9 & 438.318 & 277.8 & 440.115 \\
500 & 90  & 474.916 & 4241.8 & 474.546 & 659.0 & 475.730 \\
500 & 100 & 537.467 & 4286.3 & 536.763 & 640.7 & 538.370 \\
500 & 110 & 599.585 & 4645.3 & 599.042 & 638.6 & 600.599 \\
\hline
\end{tabular}%
}
\caption{
Comparing the Asymmetric Case: the first computation time refers to the time
in seconds for finding the estimates $\hat{f}_{0,n},\,n=1,\dots,N$,
and calculating $V_{0,0}^{test}$, when using the $\Delta$-algorithm,
while the second computation time refers to the time
in seconds for finding the estimates ${f}^\star_{0,n},\,n=1,\dots,N$, 
and calculating $V_{0,0}^\star$, when using Neural Networks with $K_{\rm train}=100000$.
}
\label{asymComparing_100000}
\end{table}

\begin{table}[ht]
\centering
\resizebox{\textwidth}{9cm}{
\begin{tabular}{ccccc}
\hline
{$D$} &                                                                            
{${\rm x}_M$} &
{\begin{tabular}[c]{@{}c@{}}CART-trees \\ \rule{0pt}{11pt} $V_{0,0}^{test}$ \\
\rule{0pt}{11pt} $K_{\rm train}=100000$  \end{tabular}} &
{\begin{tabular}[c]{@{}c@{}} computation \\ time \end{tabular}} &                           
{\begin{tabular}[c]{@{}c@{}} Neural Networks \\ $V^\star_{0,0}$ \\
\rule{0pt}{11pt} $K_{\rm train}>24\times 10^6$  \end{tabular}} \\ \hline
2   & 90  & 8.019   & 392.5 (39.1)  & 8.054 \\
2   & 100 & 13.793  & 427.3 (42.2) & 13.874 \\
2   & 110 & 21.154  & 463.8 (46.6) & 21.319 \\
3   & 90  & 11.195  & 424.0 (41.8) & 11.260 \\
3   & 100 & 18.522  & 428.2 (44.0) & 18.672 \\
3   & 110 & 27.407  & 415.2 (43.1) & 27.550 \\
5   & 90  & 16.528  & 385.0 (39.1) & 16.605 \\
5   & 100 & 26.061  & 390.9 (39.7) & 26.105 \\
5   & 110 & 36.693  & 388.7 (39.4) & 26.722 \\
10  & 90  & 26.196  & 350.6 (35.6) & 26.151 \\
10  & 100 & 38.288  & 344.3 (34.5) & 38.201 \\
10  & 110 & 50.837  & 331.6 (33.5) & 50.722 \\
20  & 90  & 37.792  & 333.8 (33.6) & 37.625 \\
20  & 100 & 51.670  & 336.0 (33.1) & 51.479 \\
20  & 110 & 65.632  & 346.6 (34.3) & 65.412 \\
30  & 90  & 44.963  & 336.6 (32.3) & 44.723 \\
30  & 100 & 59.671  & 329.5 (32.5) & 59.408 \\
30  & 110 & 74.417  & 341.7 (33.3) & 74.134 \\
50  & 90  & 54.109  & 344.0 (32.3) & 53.857 \\
50  & 100 & 69.820  & 326.2 (30.7) & 69.550 \\
50  & 110 & 85.562  & 346.1 (32.4) & 85.211 \\
100 & 90  & 66.617  & 387.1 (32.3) & 66.297 \\
100 & 100 & 83.699  & 373.9 (31.4) & 83.317 \\
100 & 110 & 100.766 & 379.4 (32.1) & 100.323 \\
200 & 90  & 79.270  & 476.7 (32.5) & 78.906 \\
200 & 100 & 97.742  & 505.8 (34.2) & 97.293 \\
200 & 110 & 116.201 & 493.2 (34.5) & 115.734 \\
500 & 90  & 96.247  & 852.8 (40.4) & 95.877 \\
500 & 100 & 116.577 & 886.1 (42.9) & 116.146 \\
500 & 110 & 136.899 & 802.6 (41.7) & 136.420 \\
\hline
\end{tabular}
}
\caption{
Symmetric Case: computation time refers to the time
in seconds for finding the estimates $\hat{f}_{0,n},\,n=1,\dots,N$,         
and calculating $V_{0,0}^{test}$, 
when using four outstanding features for the $\Delta$-algorithm.
The time in brackets shows this time without calculating $V_{0,0}^{test}$,
that is the pure time for training.
}
\label{symCompare4features}
\end{table}

\begin{table}[ht]
\centering
\resizebox{\textwidth}{9cm}{
\begin{tabular}{ccccc}
\hline
{$D$} &      
{${\rm x}_M$} &  
{\begin{tabular}[c]{@{}c@{}}CART-trees \\ \rule{0pt}{11pt} $V_{0,0}^{test}$ \\
\rule{0pt}{11pt} $K_{train}=100000$  \end{tabular}} &   
{\begin{tabular}[c]{@{}c@{}} computation \\ time \end{tabular}} &    
{\begin{tabular}[c]{@{}c@{}} Neural Networks \\ $V^\star_{0,0}$ \\
\rule{0pt}{11pt} $K_{\rm train}>24\times 10^6$  \end{tabular}} \\ \hline
2   & 90  & 14.238  & 359.2 (39.1) & 14.300 \\
2   & 100 & 19.621  & 441.5 (44.9) & 19.757 \\
2   & 110 & 26.908  & 427.2 (43.3) & 27.105 \\
3   & 90  & 18.819  & 368.5 (38.1) & 19.052 \\
3   & 100 & 26.275  & 337.0 (35.3) & 26.642 \\
3   & 110 & 35.197  & 345.7 (35.5) & 35.802 \\
5   & 90  & 27.305  & 307.4 (32.3) & 27.609 \\
5   & 100 & 37.572  & 323.1 (33.4) & 37.940 \\
5   & 110 & 48.871  & 331.1 (35.1) & 49.415 \\
10  & 90  & 85.261  & 295.1 (30.9) & 85.785 \\
10  & 100 & 103.863 & 297.8 (30.7) & 104.514 \\
10  & 110 & 122.792 & 311.5 (32.1) & 123.535 \\
20  & 90  & 125.618 & 283.8 (28.7) & 125.842 \\
20  & 100 & 149.111 & 276.5 (28.2) & 149.477 \\
20  & 110 & 172.795 & 292.0 (29.5) & 173.197 \\
30  & 90  & 154.332 & 276.3 (27.6) & 154.401 \\
30  & 100 & 181.175 & 280.8 (28.2) & 181.352 \\
30  & 110 & 207.971 & 281.7 (28.0) & 208.219 \\
50  & 90  & 196.072 & 292.9 (27.7) & 196.093 \\
50  & 100 & 227.595 & 281.6 (26.3) & 227.557 \\
50  & 110 & 259.080 & 295.9 (27.8) & 258.910 \\
100 & 90  & 263.572 & 324.4 (27.5) & 263.226 \\
100 & 100 & 302.515 & 322.0 (27.7) & 302.020 \\
100 & 110 & 341.400 & 330.2 (28.2) & 340.899 \\
200 & 90  & 344.970 & 422.4 (30.7) & 344.424 \\
200 & 100 & 393.009 & 419.9 (29.5) & 392.052 \\
200 & 110 & 440.833 & 413.7 (29.0) & 440.115 \\
500 & 90  & 477.134 & 739.1 (38.6) & 475.730 \\
500 & 100 & 539.693 & 769.2 (37.3) & 538.370 \\
500 & 110 & 602.314 & 760.5 (37.6)  & 600.599 \\
\hline
\end{tabular}
}
\caption{
Asymmetric Case: computation time refers to the time 
in seconds for finding the estimates $\hat{f}_{0,n},\,n=1,\dots,N$,  
and calculating $V_{0,0}^{test}$,                    
when using four outstanding features for the $\Delta$-algorithm.    
The time in brackets shows this time without calculating $V_{0,0}^{test}$,
that is the pure time for training.
}
\label{asymCompare4features}
\end{table}

\clearpage
\subsection{Pricing Bermudan Max-Call Options with Barrier}\label{maxCall2}
We finally compare the $\Delta$-algorithm with the tree-based algorithm 
presented in \cite{CM2020}, where the authors treat whole stopping rules as decision trees.
To emphasise the difference between our CART-trees and their trees,
we are going to say {\it Single Tree} when the optimal stopping rule itself is modeled by a tree.

The example chosen from \cite{CM2020}
deals with the optimal stopping problem (\ref{trueValue})
associated with reward function\footnote{
The parameters $r,T,N$, and $C$, have got the same meaning as in Section \ref{theBoundary}.}
\[
u(n,{\bf x})\,=\,
\exp\{-\frac{rnT}{N+1}\}
\left(\max_{1\le d\le D}{\bf x}[d]-C\,\right)^{\!+}\!\!\times x[d+1],
\quad n\in\{0,\dots,N\},\,{\bf x}\in\bfR^{D+1},
\]
and random variables $X_n,\,n=0,\dots,N$, following a stochastic process in $\bfR^{D+1}$
whose first $D$ components are independent discretised geometric Brownian motions,
while the last component keeps track of whether the maximum of the first $D$ components
has already breached an upper barrier, $\mathfrak{B}^+>0$, or not.

Following \cite{CM2020}, our training data has been generated via
\[
{\bf x}_{M+\,n}^{(k)}[d]
\,=\,
{\rm x}_M\times\exp\{\,
(\mu-\frac{\sigma^2}{2})\frac{nT}{N}
+\sigma\sqrt{\frac{T}{N}}\sum_{n'=1}^n\varepsilon_{n'}^{(k)}[d]
\,\}
\]
and
\[
{\bf x}_{M+\,n}^{(k)}[D+1]
\,=\,
{\bf 1}_{\max_{d\le D,n'\le n}{\bf x}_{M+\,n'}^{(k)}[d]\,\le\,\mathfrak{B}^+}
\]
using i.i.d.\ standard normals
$\varepsilon_{n'}^{(k)}[d],\,n'=1,\dots,N,\,k=1,\dots,K,\,d=1,\dots,D$,
while our test data following the same dynamics
has been generated using another but independent batch of i.i.d.\ standard normals
$\tilde{\varepsilon}_{n'}^{(k)}[d],\,n'=1,\dots,N,\,k=1,\dots,\tilde{K},\,d=1,\dots,D$.
For $n=0,\dots,N$,
let $X_n$ and $\tilde{X}_n$ denote the corresponding canonical random variables.

For all numerical experiments in this subsection,
the interest rate $r=0.05$, drift $\mu=0.05$, maturity $T=3$, number of time periods $N=53$,
volatility $\sigma=0.2$, strike level $C=100$, and barrier $\mathfrak{B}^+=170$,  will be fixed.
Otherwise, we consider different choices of initials ${\rm x_M}>0$, and dimensions $D$.

Table \ref{interpret-Stopping} below presents our results.

While the calculation of $V^{test}_{0,0}$ given by (\ref{ourLowerBound}) requires finding 
$\hat{f}_{0,n},\,n=1,\dots,N$,
the authors of \cite{CM2020} estimate an optimal stopping rule
by training a Single Tree. Feeding test data 
$({\bf x}_{M},\tilde{\bf x}_{M+1}^{(k)},\dots,\tilde{\bf x}_{M+N}^{(k)}),\,k=1,\dots,\tilde{K}$,
into this Single Tree yields a valid stopping-ensemble, $(\hat{\nu}^{(k)}_0)_{k=1}^{\tilde{K}}$, 
leading to the lower bound
\begin{equation}\label{singleTreeValue}
V_{0,0}^{\star\star} \,=\, \frac{1}{\tilde{K}} \sum_{k=1}^{\tilde{K} } 
u (\, \hat{\nu}^{(k)}_0 ,\, \tilde{\bf {x}}_{M+\,\hat{\nu}^{(k)}_0 }^{(k)} \,)
\end{equation}
of the true objective $\tilde{V}_{0,0}$ with respect to the given test data
(see Section \ref{problem} for the notion of stopping-ensemble).

The authors of \cite{CM2020} realised that,
even in relatively low dimensions $4\le D\le 16$,
training their Single Tree would take extremely long 
when choosing split-components from all possible dimensions, 
and hence they started looking for features.
They figured that, in the case of max-call options with barrier,
choosing split-components from only two features,
time and the rewards themselves,
would  produce good values---see the last column of Table \ref{interpret-Stopping}.

We implemented their algorithm for `growing' the Single Tree,
\cite[Algorithm 1]{CM2020},
and trained the Single Tree using training data of size $K=20000$,
only choosing split-components from the above-named two features.
Then, using test data of size $\tilde{K}=100000$,
we calculated $V^{\star\star}_{0,0}$ according to (\ref{singleTreeValue}).

Since $K$ and $\tilde{K}$ were taken from \cite{CM2020}, 
the values of $V^{\star\star}_{0,0}$ in Table \ref{interpret-Stopping}
should reproduce the corresponding values in the last column, and they do quite well,
subject to longer computation times.
However, their calculations are about $30$ times faster than ours, 
which is the well-known difference between GPU and CPU calculations.

\begin{table}[ht]
\centering
\resizebox{\textwidth}{!}{
\begin{tabular}{ccccccc}
\hline
{$D$} &
{${\rm x}_M$} &
{\begin{tabular}[c]{@{}c@{}}CART-trees \\ \rule{0pt}{11pt} $V_{0,0}^{test}$ \\
\rule{0pt}{11pt} $K_{\rm train}=100000$   \end{tabular}} &
{\begin{tabular}[c]{@{}c@{}} computation \\ time \end{tabular}} & 
{\begin{tabular}[c]{@{}c@{}} Single Tree \\ \rule{0pt}{11pt} $V_{0,0}^{\star\star}$ \\
\rule{0pt}{11pt} $K_{\rm train}=20000$   \end{tabular}} &   
{\begin{tabular}[c]{@{}c@{}} computation \\ time \end{tabular}} &    
{\begin{tabular}[c]{@{}c@{}} reported \\ in \cite{CM2020} \end{tabular}} \\ \hline
4  & 90  & 34.744 & 103.7 & 34.332 & 296.0 & 34.300(10.8) \\
4  & 100 & 43.253 & 106.7 & 43.216 & 303.0 & 43.080(6.4)  \\
4  & 110 & 49.462 & 112.4 & 49.333 & 145.1 & 49.380(4.9)  \\
8  & 90  & 45.554 & 103.0 & 45.481 & 313.1 & 45.400(7.6)  \\
8  & 100 & 51.467 & 107.1 & 51.313 & 159.0 & 51.280(3.9)  \\
8  & 110 & 54.521 & 116.7 & 54.518 & 140.7 & 54.520(3.2)  \\
16 & 90  & 51.904 & 99.2  & 51.790 & 168.9 & 51.850(4.0)  \\
16 & 100 & 54.608 & 107.4 & 54.610 & 153.7 & 54.620(3.5)  \\
16 & 110 & 55.965 & 115.8 & 55.976 & 85.4  & 56.000(2.2)  \\ \hline
\end{tabular}
}
\vspace{-0.4cm}
\caption{
The first computation time refers to the time        
in seconds for finding the estimates $\hat{f}_{0,n},\,n=1,\dots,N$, when using the $\Delta$-algorithm,  
while the second computation time refers to the time in seconds for training the Single Tree.   
The last column shows the lower bound and the corresponding computation time (in brackets)    
as reported in \cite{CM2020}.
}
\label{interpret-Stopping}
\end{table}

\vspace{-.8cm}
Interestingly, the training takes longer in lower dimensions,
which the authors of \cite{CM2020} did not explain.
When checking tree-size, we observed that the number of nodes including leaves per tree
were $11,11,7,11,7,7,7,7,5$ with respect to
$(D,{\rm x}_M)=(4,90),(4,100),(4,110)$, $(8,90),(8,100),(8,110),(16,90),(16,100),(16,110)$.
So, the longer training times in lower dimensions relate to larger trees 
being trained, but we have no intuition why the algorithm suggested in \cite{CM2020} 
would `grow' bigger trees in lower dimensions.

Comparing the performance of CART-tree-calculations and Single-Tree-calculations,
we applied the $\Delta$-algorithm with respect to the four-dimensional process
\[
(\,u(n,X_n),\,\phi_1(X_n),\,\phi_2(X_n),\,\phi_3(X_n)\,)^T,\quad n=0,\dots,N,
\]
employing the same functions\footnote{
The maxima should be taken over the first $D$ components, only, as in (\ref{ourPhi}).}
$\phi_1,\phi_2$, and $\phi_3$, defined by (\ref{ourPhi}) in the previous section.
Otherwise, 
the algorithm was combined with bagging ($B=10$) and cross-validation,    
while constraining maximum tree-depth (=10) as well as minimum node-size (=10).

We then calculated $V^{test}_{0,0}$, plugging in the \underline{same} test data
of size $\tilde{K}=100000$ we had used for the $V^{\star\star}_{0,0}$-calculations.
The computation times for $V^{test}_{0,0}$ and $V^{\star\star}_{0,0}$, 
which are both training times, are comparable,
because all calculations were done on the same machine.
Obviously, training CART-trees is significantly faster than training the Single Tree,
despite the CART-tree training data being $5$ times larger than the Single Tree training data.
Also, except for ${\rm x}_M=100,110$ in dimension $D=16$ by a very small margin,
all values of $V^{test}_{0,0}$ are bigger than those of $V^{\star\star}_{0,0}$,
and we therefore conclude that our CART-trees act better on unseen data than the Single Tree,
again applying Corollary \ref{better estimate}.
%
%
\subsection{Significance of Removal-Procedure (R) on Page \pageref{defiDelta}}\label{appli_last}
Note that, in all three examples, training data has been generated using
$\varepsilon_{n'}^{(k)}[d],\,n'=1,\dots,N,\,k=1,\dots,K,\,d=1,\dots,D$,
which are continuous random variables. As a consequence, the cardinality
$\#\{{\rm all}\;k'\}$ used in the formulation of the removal-procedure (R)
has to be ONE, leading to 
\[
\ee\!\left[\rule{0pt}{12pt}
u(\hat{\tau}_{n+1},X_{\hat{\tau}_{n+1}})-u(n,X_n)
\,|\,X_n\right]
\,=\,
u(\hat{\tau}_{n+1},X_{\hat{\tau}_{n+1}})
-
u(n,X_n),
\quad\mbox{$\pp$-a.s.},
\]
for all $n=0,\dots,N-1,\,k=1,\dots,K$,
and applying the $\Delta$-algorithm does not require performing the removal-procedure (R).

Of course, the above identity is true because the conditional expectation
is based on the empirical measure with respect to training data,
but it is unlikely to be true with respect to the 
probability measure behind the limiting generative model
(or data generating process) as mentioned in Remark \ref{bad cond exp}.
Thus, using another method of estimation than 
explicitly calculating the conditional expectation in (\ref{needed}) 
might possibly lead to an improved performance of our $\Delta$-algorithm,
which is a topic we would leave for future research.
\section{Concluding Remarks}
Our numerical experiments show that,
when estimating optimal stopping rules within the framework
developed by the authors in \cite{BCJ2019},
replacing Neural Networks by CART-trees 
leads to a plain algorithm, called $\Delta$-algorithm, 
which performs quite well at lower computational cost,
while leaving room for interpretation of the calculated rules.

As a side-effect,
applying our algorithm to pricing high-dimensional Bermudan max-call options
returned results which seem to suggest that this example
may not be complex enough
to serve as a benchmark example for high-dimensional optimal stopping.

\clearpage




\end{document}